\documentclass[4paper,twoside]{article}

\usepackage{amsmath,amsthm,amsfonts,latexsym,amscd,amssymb,enumerate}
\newcommand{\href}[2]{#2}
\swapnumbers
\theoremstyle{plain}
\newtheorem{theorem}{Theorem}[section]
\newtheorem{lemma}[theorem]{Lemma}

\newtheorem{proposition}[theorem]{Proposition}
\newtheorem{conjecture}[theorem]{Conjecture}

\theoremstyle{definition}
\newtheorem{definition}[theorem]{Definition}
\newtheorem{example}[theorem]{Example}

\theoremstyle{remark}
\newtheorem{remark}[theorem]{Remark}

\newcommand{\reals}{\mathbb{R}}
\newcommand{\complexs}{\mathbb{C}}
\newcommand{\naturals}{\mathbb{N}}
\newcommand{\integers}{\mathbb{Z}}
\newcommand{\rationals}{\mathbb{Q}}

\DeclareMathOperator{\id}{id}
\newcommand{\boundary}[1]{\partial#1}
\newcommand{\boundedops}{\mathcal{B}}

\newcommand{\abs}[1]{\left\lvert#1\right\rvert} %absolute value

\newcommand{\tensor}{\otimes}
\newcommand{\into}{\hookrightarrow}
\newcommand{\onto}{\twoheadrightarrow}
\newcommand{\iso}{\cong}

\newcommand{\subgroup}{\leq}

\DeclareMathOperator{\im}{im}      %image
\DeclareMathOperator{\coker}{coker}      %cokernel
\DeclareMathOperator{\End}{End}    %Endomorphisms
\DeclareMathOperator{\Hom}{Hom}    %Homomorphisms

\DeclareMathOperator{\tr}{tr}

\DeclareMathOperator{\ch}{ch}  % Chern character
\DeclareMathOperator{\Td}{Td}  % Todd character
\DeclareMathOperator{\ind}{ind}
\DeclareMathOperator{\sign}{sign}
\DeclareMathOperator*{\dirlim}{\varinjlim}

\newcommand{\forget}[1]{}

\newcommand{\innerprod}[1]{\langle #1 \rangle}

{\catcode`@=11\global\let\c@equation=\c@theorem}

\allowdisplaybreaks[2]
\begin{document}
\date{}%\date{Last compiled \today; last edited  \heuteIst or later}

\title{Index Theory and the Baum-Connes
    conjecture\footnote{This paper is in final form and no version of
      it is planned to be submitted for publication}}
\author{Thomas Schick\\ Mathematisches Institut\\
  Georg-August-Universit\"at G\"ottingen}

% Penn State University and Universit{\"a}t M{\"u}nster
% FB Mathematik --- Universit{\"a}t M{\"u}nster\\
% Einsteinstr.~62 --- 48159 M{\"u}nster\\
% Germany

%\vspace*{3cm}
\maketitle

%\begin{abstract}
  These notes are based on lectures on index theory, topology, and operator
  algebras at the ``School on High Dimensional Manifold Theory'' at
  the 
  ICTP in Trieste, and at the Seminari di Geometria 2002 in
  Bologna. We describe how techniques
  coming from the theory of operator algebras, in particular
  $C^*$-algebras, can be used to study manifolds. Operator algebras are extensively
  studied in their own right. We will focus on the basic definitions
  and properties, and on their relevance to the geometry and topology
  of manifolds. The link between topology and analysis is provided by
  index theorems. Starting with the classical Atiyah-Singer index
  theorem, we will explain several index theorems in detail.

  Our point of view will be in particular, that an index lives in a
  canonical way in the K-theory of a certain $C^*$-algebra. The
  geometrical context will determine, which $C^*$-algebra to use.

  A central pillar of work in the theory of $C^*$-algebras is the
  Baum-Connes conjecture. Nevertheless, it has important direct applications to the
  topology of manifolds, it implies e.g.~the Novikov conjecture. We
  will explain the Baum-Connes conjecture  and put it
  into our context.

  Several people contributed to these notes by reading preliminary
  parts and suggesting improvements, in particular Marc Johnson, Roman
  Sauer, Marco Varisco und Guido Mislin. I am very indebted to all of
  them. This is an elaboration of the first chapter of the author's
  contribution to the proceedings of the above mentioned ``School on High
  Dimensional Manifold Theory'' 2001 at the ICTP in Trieste.

%   \textbf{Important notice}: These notes are in a preliminary
%   form. They haven't gone through proofreading as thorough as the
%   author would prefer. It is certain that they will contain
%   many misprints and typos. Moreover, the presentation at many places
%   certainly can be improved. They are made public at this stage,
%   nevertheless, since I feel that they might be useful at least for
%   some participants of the school. To improve the final version, which 
%   is going to appear in the notes of the school, I would appreciate
%   any notes and comments, either in person at the school, or via email 
%   (\href{mailto:thomas.schick@math.uni-muenster.de}{e-mail: thomas.schick@math.uni-muenster.de}).

%   The current version of these notes can be found on my
%   homepage at

% \centerline{\href{http://math.uni-muenster.de/u/schickt/publ}{http://math.uni-muenster.de/u/schickt/publ}}
%  (be aware of
%   the ``t'' after schick!).
% %\end{abstract}

\section{Index theory}
\label{sec:index}

The Atiyah-Singer index theorem is one of the great achievements of
modern mathematics. It gives a formula for the index of a differential
operator (the index is by definition the dimension of the space of its
solutions minus the
dimension of the solution space for its adjoint operator) in terms
only of topological data associated to the operator and the underlying 
space. There are many good treatments of this subject available, apart 
from the original literature (most found
in~\cite{Atiyah-Collected2}). Much more detailed than the present
notes can be, because of constraints of length and time, are
e.g.~\cite{Lawson-Michelsohn(1989),Berline-Getzler-Vergne(1992),Higson-Roe(2001)}.

\subsection{Elliptic operators and their index}

We quickly review what type of operators we are looking at. This will
also fix the notation.

\begin{definition}
  Let $M$ be a smooth manifold of dimension $m$; $E,F$ smooth
  (complex) vector bundles
  on $M$. A \emph{differential operator} (of order $d$) from $E$ to
  $F$ is a $\complexs$-linear map from the space of smooth sections
  $C^\infty(E)$ of $E$ to the space of smooth sections of $F$:
  \begin{equation*}
    D\colon C^\infty(E)\to C^{\infty}(F),
  \end{equation*}
  such that in local coordinates and with local trivializations of the 
  bundles it can be written in the form
  \begin{equation*}
    D= \sum_{\abs{\alpha}\le d} A_\alpha(x)
    \frac{\partial^{\abs{\alpha}}}{\partial x^\alpha}.
  \end{equation*}
  Here $A_\alpha(x)$ is a matrix of smooth complex valued functions,
  $\alpha=(\alpha_1,\dots,\alpha_m)$ is an $m$-tuple of non-negative
  integers and
  $\abs{\alpha}=\alpha_1+\dots+\alpha_m$.
  $\partial^{\abs{\alpha}}/\partial x^\alpha$ is an abbreviation for
  $\partial^{\abs{\alpha}}/\partial x_1^{\alpha_1}\cdots\partial
  x_m^{\alpha_m}$. We require that 
  $A_\alpha(x)\ne 0$ for some $\alpha$ with $\abs{\alpha}=d$ (else,
  the operator is of order strictly smaller than $d$).

%  A \emph{real differential operator} is defined similarly, with
%  $\complexs$ replaced by $\reals$ everywhere.

Let $\pi\colon T^*M\to M$ be the bundle projection of the cotangent
bundle of $M$. We get pull-backs $\pi^*E$ and $\pi^*F$ of the bundles
$E$ and $F$, respectively, to $T^*M$.

The \emph{symbol} $\sigma(D)$ of the differential operator $D$ is the
section of the bundle
$\Hom(\pi^*E,\pi^* F)$ on $T^*M$ defined as follows:

In the above local coordinates, using $\xi=(\xi_1,\dots,\xi_m)$ as coordinate for the
cotangent vectors in $T^*M$,
in the fiber of $(x,\xi)$, the symbol $\sigma(D)$ is given by
multiplication with
\begin{equation*}
  \sum_{\abs{\alpha}=m} A_{\alpha}(x) \xi^\alpha.
\end{equation*}
Here $\xi^\alpha=\xi_1^{\alpha_1}\cdots \xi_m^{\alpha_m}$.

The operator $D$ is called \emph{elliptic}, if
$\sigma(D)_{(x,\xi)}\colon \pi^*E_{(x,\xi)}\to \pi^*F_{(x,\xi)}$ is
invertible outside the zero section of $T^*M$, i.e.~in each fiber over 
$(x,\xi)\in T^*M$ with $\xi\ne 0$. Observe that elliptic operators can 
only exist if the fiber dimensions of $E$ and $F$ coincide.

In other words, the symbol of an elliptic operator gives us two vector 
bundles over $T^*M$,  namely $\pi^*E$ and $\pi^*F$, together with a
choice of an isomorphism of the fibers of these two bundles outside
the zero section. If $M$ is compact, this gives an element of the
relative $K$-theory group $K^0(DT^*M,ST^*M)$, where $DT^*M$ and
$ST^*M$ are the disc bundle and sphere bundle of $T^*M$, respectively
(with respect to some arbitrary Riemannian metric).
\end{definition}

Recall the following definition:
\begin{definition}
  Let $X$ be a compact topological space. We define the $K$-theory of
  $X$, $K^0(X)$, to be the Grothendieck group of (isomorphism classes
  of) complex vector bundles over $X$ (with finite fiber
  dimension). More precisely, $K^0(X)$ consists of equivalence classes 
  of pairs $(E,F)$ of
  (isomorphism classes of) vector bundles over $X$, where $(E,F)\sim
  (E',F')$ if and only if there exists another vector bundle $G$ on
  $X$ such that $E\oplus F'\oplus G\iso E'\oplus F\oplus G$. One often 
  writes $[E]-[F]$ for the element of $K^0(X)$ represented by $(E,F)$.

  Let $Y$ now be a closed subspace of $X$. The \emph{relative
    $K$-theory} $K^0(X,Y)$ is given by equivalence classes of triples
  $(E,F,\phi)$, where $E$ and $F$ are complex vector bundles over $X$, 
  and $\phi\colon E|_Y\to F|_Y$ is a given isomorphism between the
  restrictions of $E$ and $F$ to $Y$. Then $(E,F,\phi)$ is isomorphic to
  $(E',F',\phi')$ if we find isomorphisms $\alpha\colon E\to E'$ and
  $\beta\colon F\to F'$ such that the following diagram commutes.
  \begin{equation*}
    \begin{CD}
      E|_Y @>{\phi}>> F|_Y\\
      @VV{\alpha}V @VV{\beta}V\\
      E'|_Y @>{\phi'}>> F'|_{Y}
    \end{CD}
  \end{equation*}
  Two pairs $(E,F,\phi)$ and $(E',F',\phi')$ are equivalent, if there
  is a bundle $G$ on $X$ such that $(E\oplus G,F\oplus
  G,\phi\oplus\id)$ is isomorphic to $(E'\oplus G,F'\oplus
  G,\phi'\oplus \id)$. 
\end{definition}

\begin{example}
  The element of $K^0(DT^*M,ST^*M)$ given by the symbol of an elliptic 
  differential operator $D$ mentioned above is represented by the
  restriction of the bundles $\pi^*E$ and $\pi^*F$ to the disc bundle
  $DT^*M$, together with the isomorphism $\sigma(D)_{(x,\xi)}\colon
  E_{(x,\xi)}\to F_{(x,\xi)}$ for $(x,\xi)\in ST^*M$.
\end{example}

\begin{example}
  Let $M=\reals^m$ and $D=\sum_{i=1}^m (\partial/\partial_i)^2$ be the 
  Laplace operator on functions. This is an elliptic differential
  operator, with symbol $\sigma(D)=\sum_{i=1}^m \xi_i^2$.

  More generally, a second-order differential operator $D\colon C^\infty(E)\to
  C^\infty(E)$ on a Riemannian manifold $M$ is a \emph{generalized
    Laplacian}, if
  $\sigma(D)_{(x,\xi)} = \abs{\xi}^2\cdot\id_{E_x}$ (the norm of the
  cotangent vector $\abs{\xi}$ is
  given by the Riemannian metric).

  Notice that all generalized Laplacians are elliptic.
\end{example}

\begin{definition}
  \emph{(Adjoint operator)}\\
  Assume that we have a differential operator $D\colon C^\infty(E)\to
  C^\infty(F)$ between two Hermitian bundles $E$ and $F$ on a
  Riemannian manifold $(M,g)$. We define an $L^2$-inner product on
  $C^\infty(E)$ by the formula 
  \begin{equation*}
\innerprod{f,g}_{L^2(E)} := \int_M \innerprod{f(x),g(x)}_{E_x}
\;d\mu(x)\qquad\forall f,g\in C^\infty_0(E),
\end{equation*}
  where $\innerprod{\cdot,\cdot}_{E_x}$ is the fiber-wise inner product 
  given by the Hermitian metric, and $d\mu$ is the measure on $M$
  induced from the Riemannian metric. Here $C^\infty_0$ is the space of
  smooth section with compact support.
  The Hilbert space completion of $C^\infty_0(E)$ with respect to this 
  inner product is called $L^2(E)$.

   The \emph{formal adjoint} $D^*$ of $D$ is then defined by
   \begin{equation*}
     \innerprod{Df,g}_{L^2(F)} =
     \innerprod{f,D^*g}_{L^2(E)}\qquad\forall f\in C^\infty_0(E),\;
       g\in C^\infty_0(F).
     \end{equation*}
   It turns out that exactly one operator with this property exists,
   which is another differential operator, and which is elliptic if
   and only if $D$ is elliptic.
\end{definition}

\begin{remark}
  The class of differential operators is quite restricted. Many
    constructions one would like to carry out with differential
    operators automatically lead out of this class. Therefore, one
    often has to use \emph{pseudodifferential operators}. 
    Pseudodifferential operators are defined as a generalization of
  differential operators. There are many well written sources dealing with the theory of
pseudodifferential operators. Since we will not discuss them in detail 
here, we omit even their precise definition and refer e.g.~to
\cite{Lawson-Michelsohn(1989)} and \cite{Shubin(1987)}.
  What we have done so far with elliptic
  operators can all be extended to pseudodifferential operators. In
  particular, they have a symbol, and the concept of ellipticity is
  defined for them. When studying elliptic differential operators,
  pseudodifferential operators naturally appear and play a very
  important role. An pseudodifferential operator $P$  (which could
    e.g.~be a differential operator) is elliptic if
  and only if a pseudodifferential operator $Q$ exists such that
  $PQ-\id$ and $QP-\id$ are so called \emph{smoothing} operators, a
  particularly nice class of pseudodifferential operators. For many
  purposes, $Q$ can be considered to act like an inverse of $P$, and
  this kind of invertibility is frequently used in the theory of elliptic
  operators. However, if $P$ happens to be an elliptic differential
  operator of positive order, then $Q$ necessarily is not a
  differential operator, but only a pseudodifferential operator.

It should be noted that almost all of the results we present here for
differential operators hold also for pseudodifferential operators, and 
often the proof is best given using them.
\end{remark}

We now want to state several important properties of elliptic
operators.

\begin{theorem}
  Let $M$ be a smooth manifold, $E$ and $F$ smooth finite dimensional
  vector bundles over $M$.
  Let $P\colon C^\infty(E)\to C^\infty(F)$ be an elliptic operator.

  Then the following holds.
  \begin{enumerate}
  \item Elliptic regularity:\\
    If $f\in L^2(E)$ is weakly in the null space of $P$,
    i.e.~$\innerprod{f,P^*g}_{L^2(E)}=0$ for all $g\in C^\infty_0(F)$, 
    then $f\in C^\infty(E)$.
  \item Decomposition into finite dimensional eigenspaces:\\
    Assume $M$ is compact and $P=P^*$ (in particular, $E=F$). Then the 
    set $s(P)$ of eigenvalues of $P$ ($P$ acting on $C^\infty(E)$) is 
    a discrete subset of $\reals$, each eigenspace $e_\lambda$
    ($\lambda\in s(P)$) is finite dimensional, and
    $L^2(E)=\oplus_{\lambda\in s(P)} e_\lambda$ (here we use the
    completed direct sum in the sense of Hilbert spaces, which means
    by definition that the
    algebraic direct sum is dense in $L^2(E)$).
  \item If $M$ is compact, then $\ker(P)$ and $\ker(P^*)$ are finite
    dimensional, and then we define the \emph{index of $P$}
    \begin{equation*}
      \ind(P):=\dim_\complexs \ker(P) - \dim_\complexs \ker(P^*).
    \end{equation*}
  \end{enumerate}
  (Here, we could replace $\ker(P^*)$ by $\coker(P)$, because these
  two vector spaces are isomorphic).
\end{theorem}

\subsection{Characteristic classes}\label{sec:char-class}

For explicit formulas for the index of a differential operator, we
will have to use characteristic classes of certain bundles
involved. Therefore, we quickly review the basics about the theory of
characteristic classes.

\begin{theorem}\label{theo:classifying_space}
  Given a compact manifold $M$ (or actually any finite CW-complex),
  there is a bijection between the isomorphism classes of
  $n$-dimensional complex vector bundles on $M$, and the set of
  homotopy classes of maps from $M$ to $BU(n)$, the classifying space
  for $n$-dimensional vector bundles. $BU(n)$ is by definition the
  space of $n$-dimensional subspaces of $\complexs^\infty$ (with an
  appropriate limit topology).

  The isomorphism is given as follows: On $BU(n)$ there is the
  tautological $n$-plane bundle $E(n)$, the fiber at each point of
  $BU(n)$ just being the subspace of $\complexs^\infty$ which
  represents this point. Any map $f\colon M\to BU(n)$ gives rise to
  the pull back bundle $f^* E(n)$ on $M$. The theorem states that each
  bundle on $M$ is isomorphic to such a pull back, and that two pull
  backs are isomorphic if and only the maps are homotopic.
\end{theorem}

\begin{definition}
  A characteristic class $c$ of vector bundles assigns to each vector
  bundle $E$ over $M$ an element $c(E)\in H^*(M)$ which is natural,
  i.e.~which satisfies
  \begin{equation*}
    c(f^*E) = f^* c(E)\qquad\forall f\colon M\to N, \quad E \text{
      vector bundle over $N$}.
  \end{equation*}
  It follows that characteristic classes are given by cohomology
  classes of $BU(n)$. 
\end{definition}

\begin{theorem}
  The integral cohomology ring  $H^*(BU(n))$ is a polynomial ring in generators
  $c_0\in H^0(BU(n))$, $c_1\in H^2(BU(n))$, $\ldots$, $c_n\in
  H^{2n}(BU(n))$. We call these generators the \emph{Chern classes} of
  the tautological bundle $E(n)$ of Theorem
  \ref{theo:classifying_space}, $c_i(E(n)):= c_i$.
\end{theorem}

\begin{definition}
  Write a complex vector bundle $E$ over $M$ as $f^*E(n)$ for $f\colon
  M\to BU(n)$ appropriate. Define
  $c_i(fE):= f^*(c_i)\in H^{2i}(M;\integers)$, this is called the
  \emph{$i$-th Chern class} of the bundle $E=f^*E(n)$.
  
  If $F$ is a real vector bundle over $M$, define the Pontryagin
  classes
  \begin{equation*}
    p_i(F):= c_{2i}(F\tensor \complexs) \in H^{4i}(M;\integers).
  \end{equation*}
  (The odd Chern classes of the complexification of a real vector
  bundle are two torsion and therefore are usually ignored).
\end{definition}

\subsubsection{Splitting principle}

\begin{theorem}\label{theo:splitting_principle}
  Given a manifold $M$ and a vector bundle $E$ over $M$, there is
  another manifold $N$ together with a map $\phi\colon N\to M$, which
  induces a monomorphism $\phi^U\colon H^*(M;\integers)\to
  H^*(N;\integers)$, and such that $\phi^*E = L^1\oplus \dots L^n$ is
  a direct sum of line bundles.
\end{theorem}

  Using Theorem \ref{theo:splitting_principle}, every question about characteristic classes of vector
  bundles can be reduced to the corresponding question for line
  bundles, and questions about the behavior under direct sums.

In particular, the following definitions makes sense:

\begin{definition}
  The Chern character is an inhomogeneous characteristic class, assigning to each complex vector
  bundle $E$ over a space $M$ a cohomology class $\ch(E)\in
  H^*(M;\rationals)$. It is characterized by the following properties:
  \begin{enumerate}
  \item Normalization: If $L$ is a complex line bundle with first
    Chern class $x$, then
    \begin{equation*}
      \ch(L) = \exp(x) = \sum_{n=0}^\infty \frac{x^n}{n!} \in H^*(M;\rationals).
    \end{equation*}
    Observe that in particular $\ch(\complexs)=1$.
  \item Additivity: $L(E\oplus F) = L(E) + L(F)$.
  \end{enumerate}
\end{definition}

\begin{proposition}
  The Chern character is not only additive, but also
  multiplicative in the following sense: for two vector bundles $E$,
  $F$ over $M$ we have
  \begin{equation*}
    \ch(E\tensor F) = \ch(E)\cup \ch(F).
  \end{equation*}
\end{proposition}

\begin{definition}
  The \emph{Hirzebruch L-class} as normalized by Atiyah and Singer is
  an inhomogeneous characteristic class, assigning to each complex vector
  bundle $E$ over a space $M$ a cohomology class $L(E)\in
  H^*(M;\rationals)$. It is characterized by the following properties:
  \begin{enumerate}
  \item Normalization: If $L$ is a complex line bundle with first
    Chern class $x$, then
    \begin{equation*}
      L(L) = \frac{x/2}{\tanh (x/2)} = 1 +\frac{1}{12} x^2 -
      \frac{1}{720} x^4 +\cdots \in H^*(M;\rationals).
    \end{equation*}
    Observe that in particular $L(\complexs)=1$.
  \item Multiplicativity: $L(E\oplus F) = L(E)L(F)$.
  \end{enumerate}
\end{definition}

\begin{definition}
  The \emph{Todd-class}  is
  an inhomogeneous characteristic class, assigning to each complex vector
  bundle $E$ over a space $M$ a cohomology class $\Td(E)\in
  H^*(M;\rationals)$. It is characterized by the following properties:
  \begin{enumerate}
  \item Normalization: If $L$ is a complex line bundle with first
    Chern class $x$, then
    \begin{equation*}
      \Td(L) = \frac{x}{1-\exp(-x)} \in H^*(X;\rationals).
    \end{equation*}
    Observe that in particular $\Td(\complexs)=1$.
  \item Multiplicativity: $L(E\oplus F) = L(E)L(F)$.
  \end{enumerate}
\end{definition}

Note that $\ch$ as well as $L$ and $\Td$ take values in the even dimensional
part 
\begin{equation*}
  H^{ev}(M;\rationals) := \oplus_{k=0}^\infty H^{2k}(M;\rationals).
\end{equation*}

\subsubsection{Chern-Weyl theory}

Chern-Weyl theory can be used to explicitly compute characteristic classes of finite dimensional
vector spaces. For a short description compare \cite{Milnor-Stasheff(1974)}. To carry out the
Chern-Weyl procedure, one has to choose a connection on the given
vector bundle $E$. This connection has a curvature $\Omega$, which is
a two form with values in the endomorphism bundle of the given vector
bundle. 

There are well defined homomorphisms
\begin{equation*}
  \sigma_r\colon \Omega^2(M;\End(E)) \to \Omega^{2r}(M;\complexs),
\end{equation*}
which can be computed in local coordinates.

\begin{theorem}
  For any finite dimension vector bundle $E$ (over a smooth manifold $M$) with connection with
  curvature $\Omega$, for the image of the $k$-th Chern class $c_k(E)$
  in cohomology with complex coefficients, we have
  \begin{equation*}
    c_k(E) = \frac{1}{(2\pi i)^k} \sigma_k(\Omega) \in H^{2k}(M;\complexs).
  \end{equation*}
  Since all other characteristic classes of complex vector bundles are
  given in terms of the Chern classes, this gives an explicit way to
  calculate arbitrary characteristic classes.
\end{theorem}

\subsubsection{Stable characteristic classes and K-theory}

The elements of $K^0(X)$ are represented by vector bundles. Therefore,
it makes sense to ask whether a characteristic class of vector bundles
can be used to define maps from $K^0(X)$ to $H^*(X)$.

It turns out, that this is not always the case. The obstacle is, that
two vector bundles $E$, $F$ represent the \emph{same} element in
$K^0(X)$ if (and only if) there is $N\in\naturals$ such that $E\oplus
\complexs^N\iso F\oplus \complexs^N$. Therefore, we have to make sure
that $c(E)=c(F)$ in this case. A characteristic class which satisfies
this property is called \emph{stable}, and evidently induces a map
\begin{equation*}
  c\colon K^0(X)\to H^*(X).
\end{equation*}

We deliberately did not specify the coefficients to be taken for
cohomology, because most stable characteristic classes will take
values in $H^*(X;\rationals)$ instead of $H^*(X;\integers)$.

The following proposition is an immediate consequence of the definition:
\begin{proposition}
  Assume a characteristic class $c$ is multiplicative, i.e.~$c(E\oplus
  F) = c(E)\cup c(F) \in H^*(X)$, and $c(\complexs)=1$. Then $c$ is a
  stable characteristic class.

  Assume a characteristic class $c$ is additive, i.e.~$c(E\oplus
  F)=c(E)+c(F)$. Then $c$ is a stable characteristic class.
\end{proposition}

It follows in particular that the Chern character, as well as Hirzebruch's $L$-class are stable characteristic classes,
i.e.~they 
define maps from the K-theory $K^0(X)\to H^*(X;\rationals)$.

The relevance of the Chern character becomes apparent by the following
theorem.

\begin{theorem}
  For a finite CW complex $X$,
  \begin{equation*}
    \ch\tensor \id_\rationals\colon K^0(X)\tensor\rationals \to
    H^{ev}(X;\rationals)\tensor\rationals = H^{ev}(X;\rationals)
  \end{equation*}
  is an isomorphism.
\end{theorem}

We have constructed \emph{relative} K-theory $K^0(X,A)$ in terms of
pairs of vector bundles on $X$ with a given isomorphism of the
restrictions to $A$. We can always find representatives such that one
of the bundles is trivialized, and the other one $E$ has in particular a
trivialization $E|_A=\complexs^n$ of its restriction to $A$. Such vector bundles
correspond to homotopy classes $[(X,A); (BU(n),pt)]$ of maps from $X$
to $BU(n)$ which map $A$ to a fixed point $pt$ in $BU(n)$.

For $k>0$, we define relative Chern classes $c_k(E, E|_A=\complexs^n)\in
H^{2k}(X,A;\integers)$ as pull back of $c_k\in H^{2k}(BU(n),pt)\iso
  H^{2k}(BU(n))$. The splitting principle also holds for such relative
  vector bundles, and therefore all the definitions we have made above
  go through in this relative situation. In particular, we can define
  a Chern character
  \begin{equation*}
    \ch\colon K^0(X,A) \to H^{ev}(X,A;\rationals).
  \end{equation*}

Given an elliptic differential operator $D$, we can apply this to our symbol element 
\begin{equation*}
  \sigma(D)\in K^0(DT^*M, ST^*M),
\end{equation*}
to obtain $\ch(\sigma(D))$.

\begin{proposition}
  Given a smooth manifold $M$ of dimension $m$, there is a
  homomorphism
  \begin{equation*}
    \pi_!\colon H^{k+m}(DT^*M,ST^*M;\reals)\to H^k(M),
  \end{equation*}
  called \emph{integration along the fiber}. It is defined as follows:
  let $\omega\in \Omega^{k+m}(DT^*M)$ be a closed differential form representing an
  element in $H^{k+m}(DT^*M,ST^*M)$ (i.e.~with vanishing restriction
  to the boundary). Locally, one can write $\omega= \sum \alpha_i\cup
  \beta_i$, where $\beta_i$ are differential forms on $M$ pulled back
  to $DT^*M$ via the projection map $\pi\colon DT^*M\to M$, and
  $\alpha_i$ are pulled back from the fiber in a local
  trivialization. Then $\pi_!\omega$ is represented by
  \begin{equation*}
    \sum_i (\int_{DT_x^*M} \alpha) \beta_i.
  \end{equation*}
  For more details about integration along the fiber, consult
  \cite[Section 6]{Bott-Tu(1982)}
\end{proposition}

\subsection{Statement of the Atiyah-Singer index theorem}

There are different variants of the Atiyah-Singer index theorem. We
start with a cohomological formula for the index.

\begin{theorem}
  Let $M$ be a compact oriented manifold of dimension $m$, and $D\colon C^\infty(E)\to
  C^{\infty}(F)$ an elliptic operator with symbol $\sigma(D)$. Define
  the \emph{Todd character} $\Td(M):=\Td(TM\tensor\complexs)\in
  H^*(M;\rationals)$. Then
  \begin{equation*}
    \ind(D) = (-1)^{m(m+1)/2} \innerprod{\pi_!\ch(\sigma(D)) \cup \Td(M), [M]} .
  \end{equation*}
  The class $[M]\in H_m(M;\rationals)$ is the fundamental class of the 
  oriented manifold $M$, and $\innerprod{\cdot,\cdot}$ is the usual
  pairing between homology and cohomology. For the characteristic
  classes, compare Subsection \ref{sec:char-class}.
\end{theorem}

If we start with specific operators
given by the geometry, explicit calculation usually give more familiar
terms on the right hand
side.

 For example, for the signature operator we obtain Hirzebruch's
signature formula expressing the signature in terms of the $L$-class,
for the Euler characteristic operator we obtain the
Gauss-Bonnet formula expressing the Euler characteristic in terms of
the Pfaffian, and for the spin or spin$^c$ Dirac operator we obtain an 
$\hat{A}$-formula. For applications, these formulas prove to be
particularly useful.

We give some more details about the signature operator, which we are
going to use later again. To define the signature operator, fix a
Riemannian metric $g$ on $M$. Assume $\dim M=4k$ is divisible by four.

The signature operator maps from a certain subspace $\Omega^+$ of the space of
differential forms to another subspace $\Omega^-$. These subspaces are 
defined as follows. Define, on $p$-forms, the operator $\tau:=
i^{p(p-1)+2k}*$, where $*$ is the Hodge-$*$ operator given by the
Riemannian metric, and $i^2=-1$. Since $\dim M$ is divisible by $4$, an
easy calculation
shows that $\tau^2=\id$. We then define $\Omega^{\pm}$ to be the
$\pm1$ eigenspaces of $\tau$.

The signature operator $D_{sig}$ is now simply defined to by
$D_{sig}:= d+d^*$, where $d$ is the exterior derivative on
differential forms, and $d^*=\pm *d*$ is its formal adjoint. We restrict
this operator to $\Omega^+$, and another easy calculation shows that
$\Omega^+$ is mapped to $\Omega^-$. $D_{sig}$ is elliptic, and a
classical calculation shows that its index is the signature of $M$
given by the intersection form in middle homology.

The Atiyah-Singer index theorem now specializes to
\begin{equation*}
  \sign(M) = \ind(D_{sig}) = \innerprod{ 2^{2k} L(TM),[M]},
\end{equation*}
with $\dim M=4k$ as above.

\begin{remark}
  One direction to generalize the Atiyah-Singer index theorem is to
  give an index formula for manifolds with boundary. Indeed, this is
  achieved in the Atiyah-Patodi-Singer index theorem. However, these
  results are much less topological than the results for manifolds
  without boundary. They are not discussed in these notes.
\end{remark}

Next, we explain the K-theoretic version of the Atiyah-Singer index
theorem. It starts with the element of $K^0(DT^*M,ST^*M)$ given by the 
symbol of an elliptic operator. Given any compact manifold $M$, there is a
well defined
homomorphism 
\begin{equation*}
K^0(DT^*M,ST^*M)\to K^0(*)=\integers,
\end{equation*}
constructed as follows. Embed $M$ into high dimensional Euclidean
space $\reals^N$. This gives an embedding of $T^*M$ into
$\reals^{2N}$, and further into its one point compatification
$S^{2N}$, with normal bundle $\nu$. In this situation, $\nu$ has
a canonical complex structure. The embedding now defines a transfer
map
\begin{equation*}
 K^0(DT^*M,ST^*M) \to K^0(S^{2N}, \infty),
\end{equation*}
by first using the Thom isomorphism to map to the (compactly
supported) K-theory of the
normal bundle, and then push forward to the K-theory of the
sphere. The latter map is given by extending a vector bundle on the
open subset $\nu$ of $S^{2N}$ which is trivialized outside a compact
set (i.e.~represents an element in compactly supported K-theory)
trivially to all of $S^{2N}$.

Compose with the Bott periodicity isomorphism to map to
$K^0(pt)=\integers$. The image of the symbol element
under this homomorphism is denoted the \emph{topological index}
$\ind_t(D)\in K^0(*)=\integers$. The reason for the terminology is
that it is obtained from the symbol only,
using purely topological constructions. The Atiyah-Singer index
theorem states that analytical and topological index coincide:
\begin{theorem}
  $ \ind_t(D) =\ind(D)$.
\end{theorem}

\subsection{The $G$-index}\label{sec:equiv-index}

\subsubsection{The representation ring}

Let $G$ be a finite group, or more generally a compact Lie group.
The representation ring $RG$ of $G$ is defined to be the Grothendieck
group of all finite dimensional complex representations of $G$,
i.e.~an element of $RG$ is a formal difference $[V]-[W]$ of two finite 
dimensional $G$-representations $V$ and $W$, and we have
$[V]-[W]=[X]-[Y]$ if and only if $V\oplus Y\iso W\oplus X$ (strictly
speaking, we have to pass to isomorphism classes of representations to 
avoid set theoretical problems). The direct sum of representations
induces the structure 
of an abelian group on $RG$, and the tensor product makes it a
commutative unital ring (the unit given by the trivial one-dimensional 
representation). 

Equivalently, we can consider $RG$ as the free abelian group generated by all isomorphism classes of finite irreducible representations (since every representation decomposes uniquely as a direct sum of irreducible ones).

To get numerical information about representations, one uses
characters: the character of a representation $\rho\colon G\to Gl(V)$
is the complex valued function
$\chi_V$ on $G$ with 
$$ \chi_V(g) = \tr (\rho(g)).$$

Elements of the representation ring can be recovered from the
corresponding characters. Therefore, equivariant index theorems are
often formulated in terms of characters.

More about this representation ring can be found
e.g.~in \cite{Broecker-Dieck(1995)}.

\subsubsection{Equivariant analytic index}

Assume now that the manifold $M$ is a compact smooth manifold with a
smooth (left)
$G$-action, and let $E,F$ be complex $G$-vector bundles on $M$. This
means that $G$ acts on $E$ and $F$ by vector bundle automorphisms
(i.e.~carries fibers to fibers linearly), and the bundle projection
maps $\pi_E\colon E\to M$ and $\pi_F\colon F\to M$ are
$G$-equivariant, i.e.~satisfy 
$$ \pi_E(ge) = g \pi_E(e)\qquad\forall g\in G,\; e\in E.$$

We assume that the actions preserve a Riemannian metric on $M$ and
Hermitian metrics on $E$ and $F$. Since $G$ is compact, this can
always be achieved by avaraging an arbitrary given metric, using a
Haar measure on $G$.

The action of $G$ on $E$ and $M$ induces actions on the spaces
$C^\infty(M)$ and $C^\infty(E)$ of smooth functions on $M$, and smooth
sections of $E$, respectively. This is given by the formulas
\begin{equation*}  
  \begin{split}
    gf (x) & = f(g^{-1}x);\qquad f\in
    C^\infty(M),\;g\in G\\
    gs(x) &= g s(g^{-1}x);\qquad s\in C^\infty(E),\; g\in G .
 \end{split}
\end{equation*}

Let $D\colon C^\infty(E)\to C^\infty(F)$ be a $G$-equivariant elliptic 
differential operator, i.e.
\begin{equation*}
  D(gs) = g D(s).
\end{equation*}
Because the action of $D$ is assumed to be isometric, the adjoint
operator $D^*$ is $G$-equivariant, as well.

Since $D$ is elliptic and $M$ is compact, the kernel and cokernel of
$D$ are finite dimensional.
If $s\in \ker(D)$, $D(gs)=gD(s)=0$, i.e.~$\ker(D)$ is a finite
dimensional $G$-representation. The same is true for $\coker(D)=\ker(D^*)$.
 We define the (analytic) $G$-index of $D$ to be
\begin{equation*}
  \ind^G(D):= [\ker(D)] - [\coker(D)] \in RG.
\end{equation*}

If $G$ is the
trivial group then $RG \iso \integers$, given by the dimension, and then
$\ind^G(D)$ evidently coincides with the usual index of
$D$.

\subsubsection{Equivariant K-theory}

Generalizing the construction of the representation ring, assume that
$X$ is a compact Hausdorff space with a $G$-action. 

Then the $G$-equivariant vector bundles over $X$ form an abelian
semigroup under direct sum, and we define the equivariant K-theory
group $K^0_G(X)$ as the Grothendieck group of this semigroup,
consisting of formal differences as in the non-equivariant
case. Similarly, if $A$ is a closed subspace of $X$ which is
$G$-invariant (i.e.~$ga\in A$ for all $a\in A$, $g\in G$), then we can
use the same recipe as in the non-equivariant case, but now with
equivariant vector bundles (and equivariant maps) to define
$  K^0_G(X,A)$. 

\begin{example}
  If $X=\{pt\}$ is the one point space with the (necessarily) trivial
  $G$-action, then an equivariant vector bundles is exactly a finite
  dimensional $G$-representation. It follows that
  \begin{equation*}
     K^0_G(pt) = RG.
  \end{equation*}
\end{example}

\subsubsection{Equivariant topological index and equivariant index theorem}

Now we define a \emph{topological}
equivariant index.

To do this, we proceed exactly in the same way as in the
non-equivariant case, but now
$G$-equivariantly. 

First observe that the symbol of an equivariant differential operator
\begin{equation*}
D\colon \Gamma(E)\to \Gamma(F)
\end{equation*}
is defined by two $G$-equivariant
vector bundles, namely $\pi^*E$ 
and $\pi^*F$ over $DT^*M$, together with the isomorphism (given by the
principal symbol of the operator) of the restrictions of $\pi^*E$ and
$\pi^*F$ to $ST^*M$. $G$-equivariance of the operator implies
immediately, that this isomorphism is $G$-equivariant as
well. Consequently, in this situation the symbol can be considered to
be an element
\begin{equation*}
  \sigma_G(D) \in K^0_G(DT^*M,ST^*M).
\end{equation*}

Next we choose an equivariant
embedding of $M$ into a suitable $G$-representation $V$. Here, we have
to assume that $G$ is compact to guarantee the existence of such an
embedding. Now, an equivariant version of the Thom isomorphism (to the
equivariant K-theory of the normal bundle) and push forward define a
transfer homomorphism
\begin{equation*}
  K^0_G(DT^*M,ST^*M)\to K^0_G( V_+,\infty),
\end{equation*}
where $V_+=V\cup \{\infty\}$ is the one point compactification of the
$G$-representation $V$
to a sphere (where the $G$-action is extended to $V_+$ by
$g\infty=\infty$ for all $g\in G$).

Last, we compose with the $G$-equivariant Bott periodicity isomorphism
\begin{equation*}
K^0_G(V_+,\infty) \xrightarrow{\iso} K^0_G(pt) =RG
\end{equation*}
to map to the
representation ring.

The image of the equivariant symbol element of $K^0_G(DT^*M,ST^*M)$
under this composition is the equivariant topological index
$\ind^G_{top}(D)$.  Again, the Atiyah-Singer index theorem says
\begin{theorem}\label{theo:equivariant_Atiyah_Singer}
  $\ind^G(D) = \ind^G_t(D) \in K_G^0(pt) = RG.$
\end{theorem}

\subsection{Families of operators and their index}

Another important generalization is given if we don't look at one operator on
one manifold, but a family of operators on a family of manifolds. More 
precisely, let $X$ be any compact topological space, $\pi\colon Y\to X$ a
locally trivial fiber bundle with fibers $Y_x:=\pi^{-1}(x)\iso M$
smooth compact manifolds ($x\in X$),
and structure group the diffeomorphisms of the typical fiber $M$. Let
$E,F$ be families 
of smooth vector bundles on $Y$ (i.e.~vector bundles which are
smooth for each fiber of the fibration $Y\to X$), and $C^\infty(E)$, $C^\infty(F)$ the continuous
sections which are smooth along the fibers. 
More precisely, $E$ and $F$ are smooth fiber bundles over $X$, the
typical fiber is a vector bundle over $M$, and the structure group
consists of diffeomorphisms of this vector bundle which are fiberwise smooth.

Assume that $D\colon
C^\infty(E)\to C^\infty(F)$ is a family $\{D_x\}$ of elliptic differential
operator along the fiber $Y_x\iso M$ ($x\in X$), i.e., in local
coordinates $D$ becomes 
\begin{equation*}
\sum_{\abs{\alpha}\le m}
A_\alpha(y,x)\frac{\partial^{\abs{\alpha}}}{\partial y^\alpha}
\end{equation*}
with $y\in M$ and $x\in X$ such that $A_\alpha(y,x)$ depends
continuously on $x$, and each $D_x$ is an elliptic differential
operator on $Y_x$.

If $\dim_\complexs{\ker(D_x)}$ is independent of $x\in X$, then all of
these vector spaces patch together to give a vector bundle called
$\ker(D)$ on $X$, and similarly for the (fiber-wise) adjoint
$D^*$. This then gives a $K$-theory element $[\ker(D)]-[\ker(D^*)]\in
K^0(X)$.

Unfortunately, it does sometimes happen that these dimensions
jump. However, using appropriate
perturbations or stabilizations, one can always define the K-theory element 
\begin{equation*}
  \ind(D):= [\ker(D)]-[\ker(D^*)] \in K^0(X),
\end{equation*}
the analytic index of the family of elliptic operators $D$. For
details on this and the following material, consult
e.g.~\cite[Paragraph 15]{Lawson-Michelsohn(1989)}.

We define the symbol of $D$ (or rather a family of symbols) exactly as
in the non-parametrized case. This gives now rise to an element in
$K^0(DT_v^*Y,ST_v^*Y)$, where $T^*_vY$ is the cotangent bundle along
the fibers. Note that all relevant spaces here are fiber bundles over
$X$, with typical fiber $T^*M$, $DT^*M$ or $ST^*M$, respectively.

Now we proceed with a family version of the construction of the topological
index, copying the construction in the non-family situation, and using 
\begin{itemize}
\item
a (fiberwise) embedding of $Y$ into $\reals^N\times X$
(which is compatible with the projection maps to $X$)
\item
the Thom isomorphism for families of vector bundles
\item the family version of Bott
periodicity,
namely
\begin{equation*}
  K^0(S^{2N}\times X,
  \{\infty\}\times X)
  \xrightarrow{\iso} K^0(X).
\end{equation*}
(Instead, one could also use the K{\"u}nneth theorem
together with ordinary Bott periodicity.)
\end{itemize}

This gives rise to $\ind_t(D)\in K^0(X)$. The Atiyah-Singer index theorem
for families states:
\begin{theorem}
  $\ind(D)=\ind_t(D)\in K^0(X)$.
\end{theorem}

The upshot of the discussion of this and the last section (for the
details the reader is referred to the literature) is that the
natural receptacle for the index of differential operators in various
situations are appropriate K-theory groups, and much of todays index
theory deals with investigating these K-theory groups.

\section{Survey on $C^*$-algebras and their $K$-theory}

More detailed references for this section are, among others,
\cite{Wegge-Olsen(1993)}, \cite{Higson-Roe(2001)}, and \cite{Blackadar(1998)}.

\subsection{$C^*$-algebras}

\begin{definition}
  A \emph{Banach algebra} $A$ is a complex algebra which is a complete normed
  space, and such that $\abs{ab}\le \abs{a}\abs{b}$ for each $a,b\in A$.

  A \emph{$*$-algebra} $A$ is a complex algebra with an anti-linear
  involution $*\colon A\to
  A$ (i.e.~$(\lambda a)^*=\overline\lambda a^*$, $(ab)^*=b^*a^*$, and
  $(a^*)^*=a$ for
  all $a,b\in A$).

  A \emph{Banach $*$-algebra} $A$ is a Banach algebra which is a
  $*$-algebra such that   $\abs{a^*}=\abs{a}$ for
  all $a\in A$. 

  A \emph{$C^*$-algebra} $A$ is a Banach $*$-algebra
  which satisfies $\abs{a^*a}=\abs{a}^2$ for all $a\in A$.

  Alternatively, a $C^*$-algebra is a Banach $*$-algebra which is
  isometrically $*$-isomorphic to a norm-closed subalgebra of the
  algebra of bounded operators on some Hilbert space $H$ (this is the
  Gelfand-Naimark representation theorem, compare
  e.g.~\cite[1.6.2]{Higson-Roe(2001)}).

  A $C^*$-algebra $A$ is called separable if there exists a countable
  dense subset of $A$.
\end{definition}

\begin{example}\label{ex:loc_comp_Cstar_algebra}
  If $X$ is a compact topological space, then $C(X)$, the algebra of
  complex valued continuous functions on $X$, is a commutative
  $C^*$-algebra (with unit). The adjoint is given by complex
  conjugation: $f^*(x) = \overline{f(x)}$, the norm is the
  supremum-norm.

  Conversely, it is a theorem that every abelian unital $C^*$-algebra
  is isomorphic to $C(X)$ for a suitable compact topological space
  $X$ \cite[Theorem 1.3.12]{Higson-Roe(2001)}.

  Assume $X$ is locally compact, and set
  \begin{equation*}
C_0(X):=\{f\colon X\to\complexs\mid f\text{ continuous},
f(x)\xrightarrow{x\to\infty} 0\}.
\end{equation*}
Here, we say $f(x)\to 0$ for $x\to\infty$, or \emph{$f$ vanishes at
  infinity}, if for all $\epsilon >0$ there is a compact subset $K$ of 
$X$ with $\abs{f(x)}<\epsilon$ whenever $x\in X-K$.
This is again a commutative $C^*$-algebra (we use the supremum norm
on $C_0(X)$), and it is unital if and only if $X$ is compact (in this
case, $C_0(X)=C(X)$).
\end{example}

\subsection{$K_0$ of a ring}

Suppose $R$ is an arbitrary ring with $1$ (not necessarily
commutative). A module $M$ over $R$ is called finitely generated
projective, if there is another $R$-module $N$ and a number $n\ge 0$
such that
\begin{equation*}
  M\oplus N \iso R^n.
\end{equation*}
This is equivalent to the assertion that the matrix ring
$M_n(R)=End_R(R^n)$ contains an idempotent $e$, i.e.~with $e^2=e$,
such that $M$ is isomorphic to the image of $e$, i.e.~$M\iso e R^n$.

\begin{example} Description of projective modules.
  \begin{enumerate}
  \item If $R$ is a field, the finitely generated projective $R$-modules
    are exactly the finite dimensional vector spaces. (In this case,
    every module is projective).
  \item If $R=\integers$, the finitely generated projective modules
    are the free abelian groups of finite rank
  \item Assume $X$ is a compact topological space and $A=C(X)$. Then,
    by the Swan-Serre theorem \cite{Swan(1962)}, $M$ is a finitely
    generated projective $A$-module if and only if $M$ is isomorphic
    to the space $\Gamma(E)$ of continuous sections of some complex vector
    bundle $E$ over $X$.
\end{enumerate}
\end{example}

\begin{definition}\label{def:K0_of_rings}
  Let $R$ be any ring with unit. $K_0(R)$ is defined to be the
  Grothendieck group of finitely generated projective modules over
  $R$, i.e.~the group of equivalence classes $[(M,N)]$ of pairs of (isomorphism 
  classes of) finitely generated projective $R$-modules $M$, $N$,
  where $(M,N)\equiv (M',N')$ if and only if there is an $n\ge 0$ with 
  \begin{equation*}
    M\oplus N'\oplus R^n\iso M'\oplus N\oplus R^n.
  \end{equation*}
  The group composition is given by 
  \begin{equation*}
    [(M,N)] + [(M',N')] := [ (M\oplus M', N\oplus N')].
  \end{equation*}
  We can think of $(M,N)$ as the formal difference of modules $M-N$.

  Any unital ring homomorphism $f\colon R\to S$ induces a map 
  \begin{equation*}
f_*\colon
  K_0(R)\to K_0(S)\colon [M]\mapsto [S\tensor_R M],
\end{equation*}
where $S$ becomes a right $R$-module via $f$. We obtain that $K_0$ is a
covariant functor from the category of unital rings to the category of 
abelian groups.
\end{definition}

\begin{example} Calculation of $K_0$.
  \begin{itemize}
  \item 
    If $R $ is a field, then $K_0(R)\iso\integers$, the isomorphism
    given by the dimension: $\dim_R(M,N):=\dim_R(M)-\dim_R(N)$.
  \item
    $K_0(\integers)\iso\integers$, given by the rank.
  \item
    If $X$ is a compact topological space, then $K_0(C(X))\iso
    K^0(X)$, the topological K-theory given in terms of complex vector
    bundles. To each vector bundle $E$ one associates the
    $C(X)$-module $\Gamma(E)$ of continuous sections of $E$.
  \item Let $G$ be a discrete group. The group algebra $\complexs G$
    is a vector space with basis $G$, and with multiplication coming
    from the group structure, i.e.~given by
    $g\cdot h = (gh)$.

    If $G$ is a finite group, then $K_0(\complexs G)$ is the complex
    representation ring of $G$.
\end{itemize}
\end{example}

\subsection{K-Theory of $C^*$-algebras}

\begin{definition}
  Let $A$ be a unital $C^*$-algebra. Then $K_0(A)$ is defined as in
  Definition \ref{def:K0_of_rings}, i.e.~by forgetting the topology of 
  $A$.
\end{definition}

\subsubsection{K-theory for non-unital $C^*$-algebras}

When studying (the K-theory of) $C^*$-algebras, one has to understand
morphisms $f\colon A\to B$. This necessarily involves studying the
kernel of $f$, which is a closed ideal of $A$, and hence a
\emph{non-unital} $C^*$-algebra. Therefore, we proceed by defining the 
$K$-theory of $C^*$-algebras without unit.

\begin{definition}
  To any $C^*$-algebra $A$, with or without unit, we assign in a
  functorial way a new, unital $C^*$-algebra $A_+$ as follows. As
  $\complexs$-vector space, $ A_+:= A\oplus \complexs$, with product
  \begin{equation*}
    (a,\lambda)(b,\mu) := (ab+\lambda a+\mu b,
    \lambda\mu)\qquad\text{for }(a,\lambda),(b,\mu)\in A\oplus\complexs.
  \end{equation*}
  The unit is given by $(0,1)$. The star-operation is defined as
  $(a,\lambda)^*:=
  (a^*,\overline{\lambda})$, and the new norm is given by
   \begin{equation*}
     \abs{(a,\lambda)} =\sup \{\abs{a x +\lambda x}\mid x\in A \text{
     with }\abs{x}=1\}
   \end{equation*}
\end{definition}

\begin{remark}
  $A$ is a closed ideal of $A_+$, the kernel of the canonical
  projection $A_+\onto \complexs$ onto the second factor. If $A$
  itself is unital, the unit of $A$ is of course different from the
  unit of $A_+$.
\end{remark}

\begin{example}
  Assume $X$ is a locally compact space, and let $X_+:=
  X\cup\{\infty\}$ be the one-point compactification of $X$. Then 
  \begin{equation*}
    C_0(X)_+ \iso C(X_+).
  \end{equation*}
  The ideal $C_0(X)$ of $C_0(X)_+$ is identified with the ideal of
  those functions
  $f\in C(X_+)$ such that $f(\infty)=0$.
\end{example}

\begin{definition}
  For an arbitrary $C^*$-algebra $A$ (not necessarily unital) define
  \begin{equation*}
    K_0(A):= \ker( K_0(A_+)\to K_0(\complexs)).
  \end{equation*}
  Any $C^*$-algebra homomorphisms $f\colon A\to B$ (not necessarily
  unital) induces a unital homomorphism $f_+\colon A_+\to B_+$. The
  induced map 
  \begin{equation*}
(f_+)_*\colon K_0(A_+)\to K_0(B_+)
\end{equation*}
maps the kernel of
  the map $K_0(A_+)\to K_0(\complexs)$ to the kernel of $K_0(B_+)\to
  K_0(\complexs)$. This means it restricts to a map $f_*\colon
  K_0(A)\to K_0(B)$. We obtain a covariant functor from the category
  of (not necessarily unital) $C^*$-algebras to abelian groups.
\end{definition}

Of course, we need the following result.
\begin{proposition}
  If $A$ is a unital $C^*$-algebra, the new and the old definition of
  $K_0(A)$ are canonically isomorphic.
\end{proposition}

\subsubsection{Higher topological K-groups}

We also want to define higher topological K-theory groups. We have an
ad hoc definition using suspensions (this is similar to the
corresponding idea in topological K-theory of spaces). For this we
need the following.

\begin{definition}
  Let $A$ be a $C^*$-algebra. We define the cone $CA$ and the
  suspension $SA$ as follows.
  \begin{equation*}
    \begin{split}
      CA &:= \{f\colon [0,1]\to A\mid f(0)=0\} \\
      SA & := \{f\colon [0,1]\to A\mid f(0)=0=f(1)\}.
  \end{split}
\end{equation*}
These are again $C^*$-algebras, using pointwise operations and the
supremum norm.

Inductively, we define
\begin{equation*}
  S^0A := A\qquad S^n A:= S(S^{n-1}A)\quad\text{for }n\ge 1.
\end{equation*}
\end{definition}

\begin{definition} Assume $A$ is a $C^*$-algebra.
  For $n\ge 0$, define
  \begin{equation*}
    K_n(A):= K_0(S^n A).
  \end{equation*}
    These are the \emph{topological K-theory groups of $A$}. For each
    $n\ge 0$, we obtain a functor from the category of $C^*$-algebras
    to the category of abelian groups.
\end{definition}

For unital $C^*$-algebras, we can also give a more direct definition
of higher K-groups (in
particular useful for $K_1$, which is then defined in terms of
(classes of) invertible matrices). This is done as follows:

\begin{definition}
  Let $A$ be a unital $C^*$-algebra. Then $Gl_n(A)$ becomes a
  topological group, and we have continuous embeddings
  \begin{equation*}
    Gl_n(A)\into Gl_{n+1}(A)\colon X\mapsto 
    \begin{pmatrix}
      X & 0\\ 0 & 1
    \end{pmatrix}.
  \end{equation*}
  We set $Gl_\infty(A):= \lim_{n\to\infty} Gl_n(A)$, and we equip
  $Gl_\infty(A)$ with the direct limit topology.
\end{definition}

\begin{proposition}
  Let $A$ be a unital $C^*$-algebra. If $k\ge 1$, then
  \begin{equation*}
    K_k(A) = \pi_{k-1}(Gl_\infty(A)) (\iso \pi_k(B Gl_\infty(A))).
  \end{equation*}

  Observe that any unital morphism $f\colon A\to B$ of unital
  $C^*$-algebras induces a map $Gl_n(A)\to Gl_n(B)$ and therefore also 
  between $\pi_k(Gl_\infty(A))$ and $\pi_k(Gl_\infty(B))$. This map
  coincides with the previously defined induced map in topological
  $K$-theory.
\end{proposition}

\begin{remark}
  Note that the topology of the $C^*$-algebra enters the definition of 
  the higher topological K-theory of $A$, and in general the
  topological K-theory of $A$ will be vastly different from the
  algebraic K-theory of the algebra underlying $A$. For connections in 
  special cases, compare \cite{Suslin-Wodzicki(1992)}.
\end{remark}

\begin{example}
  It is well known that $Gl_n(\complexs)$ is connected for each
  $n\in\naturals$. Therefore
  \begin{equation*}
    K_1(\complexs) = \pi_0(Gl_\infty(\complexs)) = 0.
  \end{equation*}
\end{example}

A very important result about $K$-theory of $C^*$-algebras is the
following long exact sequence. A proof can be found e.g.~in
\cite[Proposition 4.5.9]{Higson-Roe(2001)}.
\begin{theorem}\label{theo:long_exact}
  Assume $I$ is a closed ideal of a $C^*$-algebra $A$. Then,
  we get a short exact sequence of $C^*$-algebras $0\to I\to A\to
  A/I\to 0$, which induces a long exact sequence in K-theory
  \begin{equation*}
    \to K_n(I)\to K_n(A)\to K_n(A/I)\to K_{n-1}(I) \to \cdots \to K_0(A/I).
  \end{equation*}
\end{theorem}

\subsection{Bott periodicity and the cyclic exact sequence}

One of the most important and remarkable results about the K-theory of 
$C^*$-algebras is Bott periodicity, which can be stated as follows.

\begin{theorem}
  Assume $A$ is a $C^*$-algebra. There is a natural isomorphism,
  called the Bott map
  \begin{equation*}
    K_0(A)\to K_0(S^2A),
  \end{equation*}
  which implies immediately that there are natural isomorphism
  \begin{equation*}
    K_n(A)\iso K_{n+2}(A)\qquad\forall n\ge 0.
  \end{equation*}
\end{theorem}

\begin{remark}
  Bott periodicity allows us to define $K_n(A)$ for each
  $n\in\integers$, or to regard the K-theory of $C^*$-algebras as a
  $\integers/2$-graded theory, i.e.~to talk of $K_n(A)$ with
  $n\in\integers/2$. This way, the long exact sequence of Theorem
  \ref{theo:long_exact} becomes a (six-term) cyclic exact sequence
  \begin{equation*}
    \begin{CD}
      K_0(I) @>>> K_0(A) @>>> K_0(A/I)\\
      @AAA && @VV{\mu_*}V\\
      K_1(A/I) @<<< K_1(A) @<<< K_1(I).
    \end{CD}
  \end{equation*}
  The connecting homomorphism $\mu_*$ is the composition of the Bott
  periodicity isomorphism and the connecting homomorphism of Theorem
  \ref{theo:long_exact}. 
\end{remark}

\begin{example}

Assume $A=C(X)$ is the space of continuous functions on a compact
Hausdorff space. Particularly interesting is the case where $X$ is the
one point space.

Then 
\begin{equation*}
  S^{2N}A= \{ f\colon S^{2N}\times A\to
\complexs\mid f(*,a)=0\quad\forall
a\in A\},
\end{equation*}
where $*$ is a base point in the two
 sphere $S^{2N}$.

To see this, use that
$C(X,C(Y))=C(X\times Y)$; and the
fact that 
\begin{equation*}
  \{f\colon [0,1]^{2N}\to A\mid
  f|_{\boundary [0,1]^{2n}}=0\} =
  \{f\colon S^{2N}\to A\mid
  f(*)=0\}.
\end{equation*}

It follows that
\begin{equation*}
  K^0(S^{2N}A) \iso K^0(S^{2N}\times A,
  \{*\}\times A)
\end{equation*}
is the relative topological K-theory
of the pair of spaces $(S^{2N}\times A,
\{*\}\times A)$. 

In particular, we recover the Bott
periodicity isomorphism
\begin{equation*}
  K^0(S^{2N},\infty) \to K^0(pt);\qquad
  K^0(S^{2N}\times X, \{\infty\} \times X) \to K^0(X)
\end{equation*}
used in the definition of the
topological index and the
topological family index, respectively.
\end{example}

\subsection{The $C^*$-algebra of a group}

Let $\Gamma$ be a discrete group. Define $l^2(\Gamma)$ to be the 
Hilbert space of square summable complex valued functions on
$\Gamma$. We can write an element $f\in l^2(\Gamma)$ as a sum
$\sum_{g\in\Gamma}\lambda_g g$ with $\lambda_g\in\complexs$ and
$\sum_{g\in\Gamma}\abs{\lambda_g}^2<\infty$.

% Similarly, we define $l^2(\Gamma;\reals)$ as the real Hilbert space of
% square summable real valued functions on $\Gamma$.

We defined the \emph{complex group algebra} (often also called the
\emph{complex
group ring}) $\complexs\Gamma$ to be the complex vector space  with
basis the elements of $\Gamma$ (this can also be considered as the
space of
complex valued functions on $\Gamma$ with finite support, and as such
is a subspace of $l^2(\Gamma)$). The product in $\complexs\Gamma$ is
induced by the multiplication in $\Gamma$, namely, if
$f=\sum_{g\in\Gamma}\lambda_g g, u=\sum_{g\in\Gamma}\mu_g
g\in\complexs\Gamma$, then
\begin{equation*}
  (\sum_{g\in\Gamma}\lambda_g g)(\sum_{g\in\Gamma}\mu_g
g) := \sum_{g,h\in\Gamma} \lambda_g\mu_h (gh) =\sum_{g\in \Gamma}
\left(\sum_{h\in\Gamma}\lambda_h\mu_{h^{-1}g} \right) g .
\end{equation*}
This is a convolution product.

We have the \emph{left regular representation} $\lambda_\Gamma$ of $\Gamma$ on
$l^2(\Gamma)$, given by
\begin{equation*}
  \lambda_\Gamma(g)\cdot(\sum_{h\in \Gamma}\lambda_h h) := \sum_{h\in\Gamma}\lambda_h gh
\end{equation*}
for $g\in\Gamma$ and $\sum_{h\in\Gamma}\lambda_h h\in l^2(\Gamma)$.

This unitary representation extends linearly to $\complexs \Gamma$.

The \emph{reduced $C^*$-algebra} $C^*_r\Gamma$ of $\Gamma$ is defined
to be the norm closure of the image $\lambda_\Gamma(\complexs\Gamma)$
in the $C^*$-algebra of bounded operators on $l^2(\Gamma)$.

\begin{remark}
  It's no surprise that there is also a \emph{maximal $C^*$-algebra}
  $C^*_{max}\Gamma$ of a
  group $\Gamma$. It is defined using not only the left regular
  representation of $\Gamma$, but simultaneously all of its
  representations. We will not make use of $C^*_{max}\Gamma$ in these
  notes, and therefore will not define it here.

  Given a topological group $G$, one can define $C^*$-algebras
  $C^*_rG$ and $C^*_{max}G$ which take the topology of $G$ into
  account. They actually play an important role in the study of the
  Baum-Connes conjecture, which can be defined for (almost arbitrary)
  topological groups, but again we will not cover this subject
  here. Instead, we will throughout stick to discrete groups.
\end{remark}

\begin{example}
  If $\Gamma$ is finite, then $C^*_r\Gamma =\complexs \Gamma$ is the
  complex group ring of $\Gamma$.

  In particular, in this case $K_0(C^*_r\Gamma)\iso R\Gamma$
  coincides with the (additive group of) the complex representation
  ring of $\Gamma$.
\end{example}

\section{The Baum-Connes conjecture}

The Baum-Connes conjecture relates an object from algebraic topology,
namely the K-homology of the classifying space of a given group
$\Gamma$, to representation theory and the world of $C^*$-algebras,
namely to the K-theory of the reduced $C^*$-algebra of $\Gamma$.

Unfortunately, the material is very technical. Because of lack of
space and time we can not go into the details (even of some of the
definitions). We recommend the sources \cite{Valette(2001)}, \cite{Valette(2000)},
\cite{Higson-Roe(2001)}, \cite{Baum-Connes-Higson(1994)},
\cite{Mislin(2001)} and
\cite{Blackadar(1998)}.

\subsection{The Baum-Connes conjecture for torsion-free groups}

\begin{definition}
  Let $X$ be any CW-complex. $K_*(X)$ is the K-homology of $X$, where
  K-homology is the homology theory dual to topological K-theory. If
  $BU$ is the spectrum of topological K-theory, and $X_+$ is $X$ with
  a disjoint basepoint added, then
  \begin{equation*}
    K_n(X):= \pi_n(X_+\wedge BU).
  \end{equation*}
\end{definition}

\begin{definition}
  Let $\Gamma$ be a discrete group. A classifying space $B\Gamma$
  for $\Gamma$ is a CW-complex with the property that
  $\pi_1(B\Gamma)\iso\Gamma$, and $\pi_k(B\Gamma)=0$ if $k\ne 1$. A
  classifying space always exists, and is unique up to homotopy
  equivalence. Its universal covering $E\Gamma$ is a contractible
  CW-complex with a free cellular $\Gamma$-action, the so called
  \emph{universal space for $\Gamma$-actions}.
\end{definition}

\begin{remark}
  In the literature about the Baum-Connes conjecture, one will often
  find the definition 
  \begin{equation*}
    RK_n(X):= \dirlim K_n(Y),
  \end{equation*}
  where the limit is taken over all finite subcomplexes $Y$ of
  $X$. Note, however, that K-homology (like any homology theory in
  algebraic topology) is compatible with direct limits, which implies
  $RK_n(X)=K_n(X)$ as defined above. The confusion comes from the fact 
  that operator algebraists often use Kasparov's bivariant KK-theory to 
  define $K_*(X)$, and this coincides with the homotopy theoretic 
  definition only if $X$ is compact.
\end{remark}

Recall that a group $\Gamma$ is called torsion-free, if $g^n=1$ for
$g\in \Gamma$ and $n>0$ implies that $g=1$.

We can now formulate the Baum-Connes conjecture for torsion-free
discrete groups.
\begin{conjecture}\label{conj:torsion_free_BC}
  Assume $\Gamma$ is a torsion-free discrete group. It is known that
  there is a particular
  homomorphism, the assembly map
  \begin{equation}\label{eq:BC}
    \overline\mu_*\colon K_*(B\Gamma) \to K_*(C^*_r \Gamma)
  \end{equation}
  (which will be defined later). The \emph{Baum-Connes conjecture} says that
  this map is an isomorphism.
\end{conjecture}

\begin{example}\label{ex:finite_group_assembly_not_iso}
  The map $\overline{\mu}_*$ of Equation \eqref{eq:BC} is also defined
  if $\Gamma$ is
  not torsion-free. However, in this situation it will in general not
  be an isomorphism. This can already be seen if $\Gamma=\integers/2$. Then 
  $C^*_r \Gamma=\complexs \Gamma\iso \complexs\oplus\complexs$ as a
  $\complexs$-algebra. Consequently, 
  \begin{equation}\label{eq:K_of_Z2}
K_0(C^*_r\Gamma)\iso
  K_0(\complexs)\oplus K_0(\complexs) \iso \integers\oplus\integers.
\end{equation}
On the other hand, 
%a model for $BG$ is given by $\reals P^\infty$, the 
%infinite real projective space, and, 
using the homological Chern character, 
\begin{equation}\label{eq:H_of_Z2}
 K_0(B\Gamma)\tensor_\integers\rationals \iso \oplus_{n=0}^\infty
 H_{2n}(B\Gamma;\rationals) \iso\rationals.
\end{equation}
(Here we use the fact that the rational homology of every finite group 
is zero in positive degrees, which follows from the fact that the
transfer homomorphism $H_k(B\Gamma;\rationals)\to H_k(\{1\};\rationals)$ is
(with rational coefficients) up to a factor $\abs{\Gamma}$ a left inverse
to the map induced from
the inclusion, and therefore is injective.)

The calculations \eqref{eq:K_of_Z2} and \eqref{eq:H_of_Z2} prevent
$\mu_0$ of \eqref{eq:BC} from being an isomorphism.
\end{example}

\subsection{The Baum-Connes conjecture in general}

To account for the problem visible in Example
\ref{ex:finite_group_assembly_not_iso} if we are dealing with groups
with torsion, one replaces the left hand side by a more complicated
gadget, the equivariant K-homology of a certain $\Gamma$-space
$E(\Gamma,fin)$, the classifying space for proper actions. We will
define all of this later. Then, the Baum-Connes conjecture says the
following.
\begin{conjecture}\label{conj:general_BC}
    Assume $\Gamma$ is a discrete group. It is known that there is a particular
  homomorphism, the assembly map
  \begin{equation}\label{eq:BC_general}
    \mu_*\colon K^{\Gamma}_*(E(\Gamma,fin)) \to K_*(C^*_r \Gamma)
  \end{equation}
  (we will define it later). The conjecture says that this map is an
  isomorphism.
\end{conjecture}

\begin{remark}
  If $\Gamma$ is torsion-free, then
  $K_*(B\Gamma)=K_*^\Gamma(E(\Gamma,fin))$, and the assembly maps
  $\overline\mu$ of
  Conjectures \ref{conj:torsion_free_BC} and $\mu$ of \ref{conj:general_BC}
  coincide (see Proposition \ref{prop:Davis_Lueck_properties}).
\end{remark}

Last, we want to mention that there is also a \emph{real version} of the 
Baum-Connes conjecture, where on the left hand side the K-homology is
replaced by KO-homology, i.e.~the homology dual to the K-theory of
real vector bundles (or an equivariant version hereof), and on the right
hand side $C^*_r\Gamma$ is replaced by the real reduced
$C^*$-algebra $C^*_{r,\reals}\Gamma$.

\subsection{Consequences of the Baum-Connes conjecture}

\subsubsection{Idempotents in $C^*_r\Gamma$}\label{sec:idempotents-c_rgamma}

The connection between the Baum-Connes conjecture and idempotents is
best shown via Atiyah's $L^2$-index theorem, which we discuss first.

Given a closed manifold $M$ with an elliptic differential operator
$D\colon C^\infty(E)\to C^\infty(F)$ between two bundles on $M$, and a 
normal covering $\tilde M\to M$ (with deck transformation group
$\Gamma$, normal means that $M=\tilde M/\Gamma$), we can lift $E$, $F$ 
and $D$ to $\tilde M$, and get an elliptic $\Gamma$-equivariant
differential operator $\tilde D\colon C^\infty(\tilde E)\to
C^\infty(\tilde F)$. If $\Gamma$ is not finite, we can not use the
equivariant index of Section \ref{sec:equiv-index}. However, because
the action is free, it is possible to define an equivariant analytic index
\begin{equation*}
  \ind_\Gamma(\tilde D)\in K_{\dim M}(C^*_r\Gamma).
\end{equation*}
This is described in Example \ref{ex:Atiyah-operator_on_covering}. 

Atiyah used a certain real valued homomorphism, the $\Gamma$-dimension 
\begin{equation*}
\dim_\Gamma\colon K_0(C^*_r\Gamma)\to\reals,
\end{equation*}
to define the
$L^2$-index of $\tilde D$ (on an even dimensional manifold):
\begin{equation*}
  L^2\text{-}\ind(\tilde D):= \dim_\Gamma(\ind_\Gamma(\tilde D)).
\end{equation*}
The $L^2$-index theorem says
\begin{equation*}
  L^2\text{-}\ind(\tilde D) = \ind(D),
\end{equation*}
in particular, it follows that the $L^2$-index is an integer.

\begin{definition}
  The $\Gamma$-dimension used above can be defined as follows: an
  element $x$ of $K_0(C^*_r\Gamma)$ is given by a (formal difference
  of) finitely generated
  projective modules over $C^*_r\Gamma$. Such a module is the image of
  a projection $p\in M_n(C^*_r\Gamma)$, i.e.~a matrix $(p_{ij})$ with
  entries in $C^*_r\Gamma$ and such that $p^2=p=p^*$. $C^*_r\Gamma$ is
  by definition a certain algebra of bounded operators on
  $l^2\Gamma$. On this algebra, we can define a trace $\tr_\Gamma$ by
  $\tr_\Gamma(a)= \innerprod{a(e),e}_{l^2\Gamma}$, where $e$ is the
  function in $l^2\Gamma$ which has value one at the unit element, and
  zero everywhere else.

  We then define 
  \begin{equation*}
    \dim_\Gamma(x=[\im(p)-\im(q)]) =\sum_{i=1}^n \tr_\Gamma(p_{ii}) -\tr_\Gamma(q_{ii}).
  \end{equation*}
\end{definition}

An alternative description of the left hand side of \eqref{eq:BC} and
\eqref{eq:BC_general} shows that, as long as $\Gamma$ is torsion-free, 
the image of $\mu_0$ coincides with the subset of $K_0(C^*_r\Gamma)$
consisting of $\ind_\Gamma(\tilde D)$, where $\tilde D$ is as above. In
particular, if $\mu_0$ is surjective (and $\Gamma$ is torsion-free),
for each $x\in K_0(C^*_r\Gamma)$ we find a differential operator $D$
such that $x=\ind_\Gamma(\tilde D)$. As a consequence,
$\dim_\Gamma(x)\in\integers$, i.e.~the range of $\dim_\Gamma$ is
contained in $\integers$. This is the statement of the so called
\emph{trace conjecture}.
\begin{conjecture}\label{conj:trace_conjecture}
  Assume $\Gamma$ is a torsion-free discrete group. Then
  \begin{equation*}
    \dim_\Gamma(K_0(C^*_r\Gamma))\subset\integers.
  \end{equation*}
\end{conjecture}

On the other hand, if $x\in
K_0(C^*_r\Gamma)$ is represented by a projection $p=p^2\in
C^*_r\Gamma$, then elementary properties of $\dim_\Gamma$ (monotonicity 
and faithfulness) imply that $0\le \dim_\Gamma(p)\le 1$, and
$\dim_\Gamma(p)\notin \{0,1\}$ if $p\ne 0,1$.

Therefore, we have the following consequence of the Baum-Connes
conjecture. If $\Gamma$ is torsion-free and the Baum-Connes map $\mu_0$ is 
surjective, then $C^*_r\Gamma$ does not contain any projection
different from $0$ or $1$.

This is the assertion of the Kadison-Kaplansky conjecture:
\begin{conjecture}
  Assume $\Gamma$ is torsion-free. Then $C^*_r\Gamma$ does not contain 
  any non-trivial projections.
\end{conjecture}

The following consequence of the Kadison-Kaplansky conjecture deserves 
to be mentioned:
\begin{proposition}
  If the Kadison-Kaplansky conjecture is true for a group $\Gamma$,
then the spectrum $s(x)$ of every self adjoint element $x\in C^*_r\Gamma$ is
connected. Recall that the spectrum is defined in the following way:
\begin{equation*}
s(x):=\{\lambda\in\complexs\mid (x-\lambda\cdot 1) \text{ not
  invertible}\}. 
\end{equation*}
\end{proposition}

If $\Gamma$ is not torsion-free, it is easy to construct non-trivial
projections, and it is clear that the range of $\ind_\Gamma$ is
not contained in $\integers$. Baum and Connes originally conjectured
that it is contained in the abelian subgroup $Fin^{-1}(\Gamma)$ of
$\rationals$ generated by $\{1/\abs{F}\mid F \text{ finite subgroup of 
  }\Gamma\}$. This conjecture is not correct, as is shown by an
example of Roy \cite{Roy(1999)}. In \cite{Lueck(2001a)}, L{\"u}ck proves 
that the Baum-Connes conjecture implies that the range of
$\dim_\Gamma$ is contained in the \emph{subring} of $\rationals$
generated by $\{1/\abs{F}\mid F \text{ finite subgroup of 
  }\Gamma\}$.

\begin{remark}
  An alternative, topological proof of the fact that the
  Baum-Connes implies the Kadison-Kaplansky conjecture is given by
  Mislin in \cite{Berrick-Chatterji-Mislin(2001)}. Their proof does
  not use Atiyah's $L^2$-index theorem.
\end{remark}

\subsubsection{Obstructions to positive scalar curvature}

The Baum-Connes conjecture implies the so called ``stable
Gromov-Lawson-Rosenberg'' conjecture. This implication is a theorem
due to Stephan
Stolz. The details of this is
discussed in the lectures of Stephan Stolz \cite{Stolz(2001)},
therefore we can be very
brief. We just mention the result.

\begin{theorem}\label{theo:stable_GLR}
  Fix a group $\Gamma$.
 Assume that $\mu$ in the real version of \eqref{eq:BC_general}
 discussed in Section \ref{sec:real-c-algebras} is
  injective (which follows e.g.~if $\mu$ in \eqref{eq:BC_general} is
  an isomorphism), and assume that $M$ is a closed spin manifold with
  $\pi_1(M)=\Gamma$. Assume that a certain (index theoretic)
  invariant $\alpha(M)\in K_{\dim M}(C^*_{\reals,r}\Gamma)$
  vanishes. Then there is an $n\ge 0$ such that $M\times B^n$ admits a 
  metric with positive scalar curvature.
\end{theorem}

Here, $B$ is any simply connected $8$-dimensional spin manifold with
$\hat{A}(M)=1$. Such a manifold is called a \emph{Bott manifold}.

The converse of Theorem \ref{theo:stable_GLR}, i.e.~positive scalar
curvature implies vanishing of $\alpha(M)$, is true for arbitrary
groups and without knowing anything about the Baum-Connes conjecture.

\subsubsection{The Novikov conjecture about higher signatures}

\paragraph{Direct approach}

The original form of the Novikov conjecture states that higher
signatures are homotopy invariant.

More precisely, let $M$ be an (even dimensional) closed oriented manifold
with fundamental group $\Gamma$. Let $B\Gamma$ be a classifying space
for $\Gamma$. There is a unique (up to homotopy) \emph{classifying
  map} $u\colon M\to B\Gamma$ which is defined by the property that it 
induces an isomorphism on $\pi_1$. Equivalently, $u$ classifies a
universal covering of $M$.

Let $L(M)\in H^*(M;\rationals)$ be the Hirzebruch L-class (as
normalized by Atiyah and Singer). Given any cohomology class $a\in
H^*(B\Gamma,\rationals)$, we define the higher signature
\begin{equation*}
  \sigma_a(M):= \innerprod{L(M)\cup u^*a, [M]} \in\rationals.
\end{equation*}
Here $[M]\in H_{\dim M}(M;\rationals)$ is the fundamental class of the 
oriented manifold $M$, and $\innerprod{\cdot,\cdot}$ is the usual
pairing between cohomology and homology.

Recall that the Hirzebruch signature theorem states that $\sigma_1(M)$
is the signature of $M$, which evidently is an oriented homotopy invariant.

The Novikov conjecture generalizes this as follows.
\begin{conjecture}
  Assume $f\colon M\to M'$ is an oriented homotopy equivalence between 
  two even dimensional closed oriented manifolds, with (common)
  fundamental group
  $\pi$. ``Oriented'' means that $f_*[M]=[M']$. Then all higher
  signatures of $M$ and
  $M'$ are equal, i.e.
  \begin{equation*}
    \sigma_a(M) = \sigma_a(M')\qquad\forall a\in H^*(B\Gamma,\rationals).
  \end{equation*}
\end{conjecture}

There is an equivalent reformulation of this conjecture in terms of K-ho\-mo\-lo\-gy.
To see this, let $D$ be the signature operator of $M$. (We assume here
that $M$ is smooth, and we choose a Riemannian metric on $M$ to define 
this operator. It is an elliptic differential operator on $M$.)
The operator $D$ defines an element in the K-homology of $M$, $[D]\in
K_{\dim M}(M)$. Using the map $u$, we can push $[D]$ to $K_{\dim
  M}(B\Gamma)$. We define the higher signature $\sigma(M):= u_*[D] \in 
K_{\dim M}(B\Gamma)\tensor\rationals$. It turns out that 
\begin{equation*}
  2^{\dim M/2} \sigma_a(M) =\innerprod{a, ch (\sigma(M))}\qquad\forall 
  a\in H^*(B\Gamma;\rationals),
\end{equation*}
where $ch\colon K_*(B\Gamma)\tensor\rationals \to
H_*(B\Gamma,\rationals)$ is the homological Chern character (an
isomorphism).

Therefore, the Novikov conjecture translates to the statement that
$\sigma(M)=\sigma(M')$ if $M$ and $M'$ are oriented homotopy equivalent.

Now one can show \emph{directly} that
\begin{equation*}
  \overline\mu(\sigma(M)) = \overline\mu(\sigma(M')) \in K_*(C^*_r\Gamma),
\end{equation*}
if $M$ and $M'$ are oriented homotopy equivalent. Consequently,
rational injectivity of the Baum-Connes map $\overline\mu$ immediately implies the
Novikov conjecture. If $\Gamma$ is torsion-free, this is part of the
assertion of the Baum-Connes conjecture. Because of this relation,
injectivity of the Baum-Connes
map $\mu$ is often called the ``analytic Novikov conjecture''.

\paragraph{Groups with torsion}

For an arbitrary group $\Gamma$, we have a factorization of $\overline\mu$ as
follows:
\begin{equation*}
  K_*(B\Gamma)\xrightarrow{f} K_*^\Gamma(E(\Gamma,fin)) \xrightarrow{\mu}
  K_*(C^*_r\Gamma).
\end{equation*}
One can show that $f$ is rationally injective, so that rational
injectivity of the Baum-Connes map $\mu$ implies the Novikov
conjecture also in general.

\subsection{The universal space for proper actions}

\begin{definition}
  Let $\Gamma$ be a discrete group and $X$ a Hausdorff space with an
  action of $\Gamma$. We say that the action is \emph{proper}, if
  for all $x,y\in X$ there are open neighborhood $U_x\ni x$ and
  $U_y\ni y$ such that $gU_x\cap U_y$ is non-empty only for finitely
  many $g\in \Gamma$ (the number depending on $x$ and $y$).

  The action is said to be \emph{cocompact}, if $X/\Gamma$ is compact.
\end{definition}

\begin{lemma}
  If the action of $\Gamma$ on $X$ is proper, then for each $x\in X$
  the \emph{isotropy group} $\Gamma_x:=\{g\in\Gamma\mid gx=x\}$ is
  finite.
\end{lemma}

\begin{definition}
    Let $\Gamma$ be a discrete group. A CW-complex $X$ is a
    \emph{$\Gamma$-CW-complex}, if $X$ is a CW-complex with a
  cellular action of $\Gamma$ with the additional property that,
  whenever $g(D)\subset D$ for a cell $D$ of $X$ and some
  $g\in\Gamma$, then $g|_D=\id_D$, i.e.~$g$ doesn't move $D$ at all.
\end{definition}

\begin{remark}
  There exists also the notion of $G$-CW-complex for topological
  groups $G$ (taking the topology of $G$ into account). These have to
  be defined in a different way, namely by gluing together
  $G$-equivariant cells $D^n\times G/H$. In general, such a
  $G$-CW-complex is not an ordinary CW-complex.
\end{remark}

\begin{lemma}
  The action of a discrete group $\Gamma$ on a $\Gamma$-CW-complex is
  proper if and only if every isotropy group is finite.
\end{lemma}

\begin{definition}\label{def:univesal_space_for_proper_actions}
  A proper $\Gamma$-CW-complex $X$ is called \emph{universal}, or more
  precisely \emph{universal for proper actions}, if for every proper
  $\Gamma$-CW-complex $Y$ there is a $\Gamma$-equivariant map $f\colon 
  Y\to X$ which is unique up to $\Gamma$-equivariant homotopy. Any
  such space is denoted $E(\Gamma,fin)$ or $\underline{E}\Gamma$.
\end{definition}

\begin{proposition}
  A $\Gamma$-CW-complex $X$ is universal for proper actions if and
  only if the fixed point set
  \begin{equation*}
X^H:=\{x\in X\mid hx=x\quad\forall h\in H\}
\end{equation*}
is empty whenever $H$ is an infinite subgroup of $\Gamma$, and is
contractible (and in particular non-empty) if $H$ is a finite subgroup 
of $\Gamma$.
\end{proposition}

\begin{proposition}
  If $\Gamma$ is a discrete group, then $E(\Gamma,fin)$ exists and is
  unique up to $\Gamma$-homotopy equivalence.
\end{proposition}

\begin{remark}
  The general context for this discussion are actions of a group
  $\Gamma$ where the isotropy belongs to a fixed family of subgroups
  of $\Gamma$ (in our case, the family of all finite subgroups). For
  more information, compare \cite{Dieck(1972)}.
\end{remark}

\begin{example}\strut
  \begin{itemize}
  \item If $\Gamma$ is torsion-free, then $E(\Gamma,fin)=E\Gamma$, the 
    universal covering of the classifying space $B\Gamma$. Indeed,
    $\Gamma$ acts freely on $E\Gamma$, and $E\Gamma$ is contractible.
  \item If $\Gamma$ is finite, then $E(\Gamma,fin)=\{*\}$.
  \item If $G$ is a connected Lie group with maximal compact subgroup
    $K$, and $\Gamma$ is a discrete subgroup of $G$, then
    $E(\Gamma,fin)= G/K$ \cite[Section 2]{Baum-Connes-Higson(1994)}.
  \end{itemize}
\end{example}

\begin{remark}
  In the literature (in particular, in
  \cite{Baum-Connes-Higson(1994)}), also a slightly different notion of
  universal
  spaces is discussed. One allows $X$ to be any proper metrizable
  $\Gamma$-space, and requires the universal property for all proper
  metrizable $\Gamma$-spaces $Y$. For discrete groups (which are the
  only groups we are discussing here), a
  universal space in the sense of Definition
  \ref{def:univesal_space_for_proper_actions} is universal in this
  sense.

  However, for some of the proofs
  of the Baum-Connes
  conjecture (for special groups) it is useful to use certain models
  of $E(\Gamma,fin)$ (in the broader sense) coming from the geometry
  of the group, which are not
  $\Gamma$-CW-complexes. 
\end{remark}

\subsection{Equivariant K-homology}

Let $\Gamma$ be a discrete group. We have seen that, if $\Gamma$ is
not torsion-free, the
assembly map \eqref{eq:BC} is not an isomorphism. To account for that, 
we replace $K_*(B\Gamma)$ by the equivariant K-theory of
$E(\Gamma,fin)$. Let $X$ be any proper 
$\Gamma$-CW complex. The original definition of equivariant K-homology is
due to Kasparov, making ideas of Atiyah precise. In this definition,
elements of $K_*^\Gamma(X)$ are equivalence classes of generalized
elliptic operators. In \cite{Davis-Lueck(1998)}, a more homotopy
theoretic definition of $K_*^\Gamma(X)$ is given, which puts the
Baum-Connes conjecture in the context of other isomorphism conjectures.

\subsubsection[Homotopy theoretic equivariant K-homology]{Homotopy
  theoretic definition of equivariant K-homology}
\label{sec:homot-ther-defin}

The details of this definition are quite technical, using spaces and
spectra over the orbit category of the discrete group $\Gamma$. The
objects of the orbit category are the orbits $\Gamma/H$, $H$ any
subgroup of $\Gamma$. The morphisms from $\Gamma/H$ to $\Gamma/K$ are
simply the $\Gamma$-equivariant maps.

In this section, spectra are used in the sense of homotopy theory,
they are a generalization of topological spaces, in particular of
CW-complexes. For a basic introduction to this theory, one may consult
e.g.~\cite[Chapter 3]{Kochman(1996)} or \cite{Davis-Kirk(2001)}. For
the (more intricate) constructions
mentioned in here, only the original literature
\cite{Davis-Lueck(1998)} is available.

In this setting, any spectrum over the orbit category gives rise to an
equivariant homology theory. The decisive step is then the
construction of a (periodic) topological K-theory spectrum
$\mathbf{K}^\Gamma$ over the orbit category of $\Gamma$. This gives
us then a functor from the category of (arbitrary)
$\Gamma$-CW-complexes to the category of (graded) abelian groups, the
\emph{equivariant $K$-homology} $K^\Gamma_*(X)$ ($X$ any
$\Gamma$-CW-complex).

The important property (which justifies the name ``topological
K-theory spectrum) is that
\begin{equation*}
  K^\Gamma_k(\Gamma/H) = \pi_k(\mathbf{K}^\Gamma(\Gamma/H)) \iso K_k(C^*_rH)
\end{equation*}
for every subgroup $H$ of $\Gamma$. In particular,
\begin{equation*}
  K^\Gamma_k(\{*\}) \iso K_k(C^*_r\Gamma).
\end{equation*}

Moreover, we have the following properties:
\begin{proposition}\label{prop:Davis_Lueck_properties}
\begin{enumerate}
\item Assume $\Gamma$ is the trivial group. Then
  \begin{equation*}
    K^\Gamma_*(X) = K_*(X),
  \end{equation*}
  i.e.~we get back the ordinary K-homology introduced above.
\item If $H\subgroup \Gamma$ and $X$ is an $H$-CW-complex, then there
  is a natural isomorphism
      \begin{equation*}
        K_*^H(X) \iso K_*^\Gamma(\Gamma\times_H X).
      \end{equation*}
      Here $\Gamma\times_H X = \Gamma\times H/\sim$, where we divide
      out the equivalence relation generated by $(gh,x)\sim (g,hx)$
      for $g\in\Gamma$, $h\in H$ and $x\in X$. This is in the obvious
      way a left $\Gamma$-space.
    \item Assume $X$ is a free $\Gamma$-CW-complex. Then there is a
      natural isomorphism
  \begin{equation*}
    K_*(\Gamma\backslash X)\to K_*^\Gamma(X).
  \end{equation*}
  In particular, using the canonical $\Gamma$-equivariant map
  $E\Gamma\to E(\Gamma,fin)$, we get a natural homomorphism
  \begin{equation*}
    K_*(B\Gamma) \xrightarrow{\iso} K_*^\Gamma(E\Gamma)\to K_*^\Gamma(E(\Gamma,fin)).
  \end{equation*}
\end{enumerate}
\end{proposition}

\subsubsection[Analytic equivariant K-homology]{Analytic definition
  of equivariant K-homology}

Here we will give the original definition, which embeds into the
powerful framework of equivariant KK-theory, and which is used for
almost all proofs of special cases of the Baum-Connes
conjecture. However, to derive some of the consequences of the
Baum-Connes conjecture, most notably about the positive scalar
curvature question ---this is discussed in one of the lectures of 
Stephan Stolz \cite{Stolz(2001)}--- the homotopy theoretic definition is used.

\begin{definition}
  A Hilbert space $H$ is called \emph{($\integers/2$)-graded}, if $H$
  comes with an orthogonal sum
  decomposition $H=H_0\oplus H_1$. Equivalently, a
  unitary operator $\epsilon$ with $\epsilon^2=1$ is given on
  $H$. The subspaces $H_0$ and $H_1$ can be recovered as the 
  $+1$ and $-1$ eigenspaces of $\epsilon$, respectively.

  A bounded operator $T\colon H\to H$ is called \emph{even} (with respect to the
  given grading), if $T$ commutes with $\epsilon$, and \emph{odd}, if
  $\epsilon$ and
  $T$ anti-commute, i.e.~if $T\epsilon=-\epsilon T$. An even operator 
  decomposes as $T=\left(
    \begin{smallmatrix}
      T_0 & 0 \\ 0 & T_1
    \end{smallmatrix}\right)$, an odd one as $T= \left(
    \begin{smallmatrix}
      0 & T_0\\ T_1 & 0
    \end{smallmatrix}\right)$ in the given decomposition $H=H_0\oplus H_1$.
\end{definition}

\begin{definition}
  A \emph{generalized elliptic $\Gamma$-operator on $X$}, or a \emph{cycle
    for $\Gamma$-K-homology of the $\Gamma$-space $X$}, simply a \emph{cycle} for short, is a triple
  $(H,\pi,F)$, where 
  \begin{itemize}
  \item $H=H_0\oplus H_1$ is a $\integers/2$-graded $\Gamma$-Hilbert space
    (i.e.~the direct sum of two Hilbert spaces with unitary
    $\Gamma$-action)
  \item $\pi$ is a $\Gamma$-equivariant $*$-representation of
    $C_0(X)$ on even bounded operators of $H$ (equivariant means that
    $\pi(fg^{-1}) = g\pi(f)g^{-1}$ for all $f\in C_0(X)$ and all
    $g\in\Gamma$.
%graded means, since $C_0(X)$ is trivially graded, that
%    $\pi(f)$ respects the grading for each $f\in C_0(X)$).
  \item $F\colon H\to H$ is a bounded, $\Gamma$-equivariant, self adjoint
    operator such that $\pi(f) (F^2-1)$ and
    $[\pi(f),F]:=\pi(f)F-F\pi(f)$ are compact operators for all $f\in
    C_0(X)$.
    Moreover, we require that $F$ is odd, i.e.~$F=
\left(    \begin{smallmatrix}
      0 & D^*\\ D & 0
    \end{smallmatrix}\right)$
    in the decomposition $H=H_0\oplus H_1$.
\end{itemize}
\end{definition}

\begin{remark}
  There are many different definitions of cycles, slightly
  weakening or strengthening some of the conditions. Of course, this
  does not effect the equivariant K-homology groups which are
  eventually defined using them.
\end{remark}

\begin{definition}
  We define the direct sum of two cycles in the obvious way.
\end{definition}

\begin{definition}
  Assume $\alpha=(H,\pi,F)$ and $\alpha'=(H',\pi',F')$ are two cycles.
  \begin{enumerate}
  \item They are called (isometrically) isomorphic, if there is a
    $\Gamma$-equivariant grading preserving isometry $\Psi\colon H\to
    H'$ such that $\Psi\circ \pi(f)= \pi'(f)\circ \Psi$ for all $f\in
    C_0(X)$ and $\Psi\circ F=F'\circ \Psi$.
  \item They are called \emph{homotopic} (or \emph{operator
      homotopic}) if $H=H'$, $\pi=\pi'$, and there is a norm
    continuous path $(F_t)_{t\in[0,1]}$ of operators with $F_0=F$ and
    $F_1=F'$ and such that $(H,\pi,F_t)$ is a cycle for each $t\in
    [0,1]$.
  \item $(H,\pi,F)$ is called \emph{degenerate}, if $[\pi(f),F]=0$ and 
    $\pi(f)(F^2-1)=0$ for each $f\in C_0(X)$.
  \item The two cycles are called equivalent if there are degenerate
    cycles $\beta$ and $\beta'$ such that $\alpha\oplus\beta$ is
    operator homotopic to a cycle isometrically isomorphic to
    $\alpha'\oplus\beta'$.
  \end{enumerate}

  The set of equivalence classes of cycles is denoted
  $KK_0^\Gamma(X)$. (Caution, this is slightly unusual, mostly one
  will find the notation $K^\Gamma(X)$ instead of $KK^\Gamma(X)$).
\end{definition}

\begin{proposition}
  Direct sum induces the structure of an abelian group on $KK_0^\Gamma(X)$.
\end{proposition}

\begin{proposition}
  Any proper $\Gamma$-equivariant map $\phi\colon X\to Y$ between two
  proper $\Gamma$-CW-complexes induces a homomorphism
  \begin{equation*}
    KK^\Gamma_0(X)\to KK^\Gamma_0(Y)
  \end{equation*}
  by $(H,\pi,F)\mapsto (H,\pi\circ \phi^*,F)$, where $\phi^*\colon
  C_0(Y)\to C_0(X)\colon f\mapsto f\circ \phi$ is defined since $\phi$
  is a proper map (else
  $f\circ \phi$ does not necessarily vanish at infinity).
\end{proposition}
Recall that a continuous map $\phi\colon X\to Y$ is called
\emph{proper} if the inverse image of every compact subset of $Y$ is
compact .

  It turns out that the analytic definition of equivariant K-homology
  is quite flexible. It is designed to make it easy to construct
  elements of these groups ---in many geometric situations they
  automatically show up. We give one of the most typical examples of
  such a situation, which we will have used in Section
  \ref{sec:idempotents-c_rgamma}..

We need the following definition:
\begin{definition}\label{def:L2_sections}
  Let $\overline M$ be a (not necessarily compact) Riemannian manifold 
  without boundary, which is complete as a metric space. Define
  \begin{equation*}
    L^2\Omega^p(\overline M) := \{\omega\text{ measurable $p$-form on
      }M\mid \int_{\overline M} \abs{\omega(x)}^2_x \;d\mu(x) <\infty\}.
  \end{equation*}
  Here, $\abs{\omega(x)}_x$ is the pointwise norm (at $x\in\overline
  M$) of $\omega(x)$, which is given by the Riemannian metric, and
  $d\mu(x)$ is the measure induced by the Riemannian metric.

  $L^2\Omega^p(\overline M)$ can be considered as the Hilbert space
  completion of the space of compactly supported $p$-forms on
  $\overline M$. The inner product is given by integrating the
  pointwise inner product, i.e.
  \begin{equation*}
    \innerprod{\omega,\eta}_{L^2} := \int_{\overline M}
    \innerprod{\omega(x),\eta(x)}_x\; d\mu(x).
  \end{equation*}
\end{definition}

\begin{example}\label{ex:Atiyah-operator_on_covering}

  Assume that
  $M$ is a compact even dimensional Riemannian manifold. Let
  $X=\overline M$ be a normal covering of $M$ with deck transformation group
   $\Gamma$ (normal means that $X/\Gamma=M$). Of course, the action is
   free, in particular, proper. Let $E=E_0\oplus E_1$ be a graded
  Hermitian vector bundle on $M$, and 
  \begin{equation*}
D\colon C^\infty(E)\to C^\infty(E)
\end{equation*}
an
  odd elliptic self adjoint differential operator (odd means that $D$
  maps the subspace
  $C^\infty(E_0)$ to $C^\infty(E_1)$, and vice versa). If $M$ is
  oriented, the signature operator on $M$ is such an operator, if $M$
  is a spin-manifold, the same is true for its Dirac operator.

  Now we can pull back $E$ to a bundle $\overline E$ on $\overline M$, and
  lift $D$ to an operator $\overline D$ on $\overline E$. The assumptions
  imply that $\overline D$ extends to an unbounded self adjoint operator
  on $L^2(\overline E)$, the space of square integrable sections of
  $\overline E$. This space is the completion of $C^\infty_c(\overline
  E)$ with
  respect to the canonical inner product (compare Definition
  \ref{def:L2_sections}). (The
  subscript c denotes sections with compact support). Using the
  functional calculus, we can replace
  $\overline D$ by 
  \begin{equation*}
F:= (\overline D^2+1)^{-1/2} \overline D\colon L^2(\overline
  E)\to L^2(\overline E).
\end{equation*}
Observe that
\begin{equation*}
 L^2(\overline E)=L^2(\overline
  E_0)\oplus L^2(\overline E_1)
\end{equation*}
is a $\integers/2$-graded Hilbert
  space with a unitary
  $\Gamma$-action, which admits an (equivariant) action $\pi$ of
  $C_0(\overline M)=C_0(X)$ by fiber-wise multiplication. This action
  preserves
  the grading. Moreover, $\overline D$ as well as $F$ are odd,
  $\Gamma$-equivariant, self adjoint operators on $L^2(\overline E)$ and $F$
  is a bounded operator. From ellipticity it follows that
  \begin{equation*}
  \pi(f)(F^2-1) = -\pi(f) 
  (\overline D^2+1)^{-1} 
\end{equation*}
is compact for each $f\in C_0(\overline M)$
  (observe that this is not true for $(\overline D^2+1)^{-1}$ itself, if
  $\overline M$ is not compact). Consequently, $(L^2(E),\pi,F)$ defines
  an (even) cycle for $\Gamma$-K-homology, i.e.~it represents an element 
  in $KK^\Gamma_0(X)$.

  One can slightly reformulate the construction as follows: $\overline M$ 
  is a principal $\Gamma$-bundle over $M$, and $l^2(\Gamma)$ has a (unitary)
  left $\Gamma$-action. We therefore can construct the associated flat 
  bundle
  \begin{equation*}
 L:=l^2(\Gamma)\times_\Gamma\overline M
\end{equation*}
on $M$ with fiber
  $l^2(\Gamma)$. Now we can twist $D$ with this bundle $L$,
  i.e.~define 
  \begin{equation*}
\overline D:=\nabla_L\tensor \id+\id\tensor D\colon
  C^\infty(L\tensor E)\to C^\infty(L\tensor E),
\end{equation*}
using the given flat
  connection $\nabla_L$ on $L$. Again, we can
  complete to $L^2(L\tensor E)$ and define 
  \begin{equation*}
F:=(\overline
  D^2+1)^{-1/2}\overline D.
\end{equation*}
The left action of $\Gamma$ on $l^2\Gamma$
  induces an action of $\Gamma$ on $L$ and then a unitary action on
  $L^2(L\tensor E)$. Since $\nabla_L$ preserves the $\Gamma$-action,
  $\overline D$ is $\Gamma$-equivariant.
  There is
  a canonical $\Gamma$-isometry
  between $L^2(L\tensor E)$ and $L^2(\overline E)$ which identifies the
  two versions of $\overline D$ and $F$. The action of $C_0(\overline M)$ on 
  $L^2(L\tensor E)$ can be described by identifying $C_0(\overline M)$ with the
  continuous sections of $M$ on the associated bundle
  \begin{equation*}
C_0(\Gamma)\times_\Gamma\overline M,
\end{equation*}
where $C_0(\Gamma)$ is the
  $C^*$-algebra of functions on $\Gamma$ vanishing at infinity, and
  then using the obvious action of $C_0(\Gamma)$ on $l^2(\Gamma)$.

  It is easy to see how this examples generalizes to
  $\Gamma$-equivariant elliptic differential operators on manifolds
  with a proper, but not necessarily free, $\Gamma$-action (with the
  exception of the last part, of course).

  Work in progress of Baum, Higson and Schick \cite{Baum-Higson-Schick(2001)}
  suggests the (somewhat surprising) fact that, given any proper
  $\Gamma$-CW-complex $Y$, we can,
  for each element $y\in KK^\Gamma_0(Y)$, find such a proper
  $\Gamma$-manifold $X$, together with a $\Gamma$-equivariant map
  $f\colon X\to Y$ and an elliptic differential operator on $X$ giving
   an element $x\in KK^\Gamma_0(X)$ as in the example, such that
  $y=f_*(x)$.
\end{example}

Analytic K-homology is homotopy invariant, a proof can be found in
\cite{Blackadar(1998)}.
\begin{theorem}
  If $\phi_1,\phi_2\colon X\to Y$ are proper $\Gamma$-equivariant maps 
  which are homotopic through proper $\Gamma$-equivariant maps,
  then
  \begin{equation*}
    (\phi_1)_*=(\phi_2)_*\colon KK_*^\Gamma(X)\to KK_*^\Gamma(Y).
  \end{equation*}
\end{theorem}

\begin{theorem}
  If $\Gamma$ acts freely on $X$, then
  \begin{equation*}
    KK_*^\Gamma(X)\iso K_*(\Gamma\backslash X),
  \end{equation*}
  where the right hand side is the ordinary K-homology of
  $\Gamma\backslash X$.
\end{theorem}

\begin{definition}
  Assume $Y$ is an arbitrary proper $\Gamma$-CW-complex. Set
  \begin{equation*}
    RK^\Gamma_*(Y):= \dirlim KK^\Gamma_*(X),
  \end{equation*}
  where we take the direct limit over the direct system of
  $\Gamma$-invariant subcomplexes of $Y$ with compact quotient (by the
  action of $\Gamma$).
\end{definition}

\begin{definition}
  To define higher (analytic) equivariant K-homology, there are two
  ways. The short one only works for complex K-homology. One
  considers cycles and an
  equivalence relation exactly as above --- with the notable exception 
  that one does not require any grading! This way, one defines
  $KK^\Gamma_1(X)$. Because of Bott periodicity (which has period $2$), 
  this is enough to define all K-homology groups ($KK_n^\Gamma(X)=
  KK_{n+2k}^\Gamma(X)$ for any $k\in\integers$).

  A perhaps more conceptual approach is the following. Here, one
  generalizes the notion of a graded Hilbert space by the notion of a
  $p$-multigraded Hilbert space ($p\ge 0$). This means that the graded 
  Hilbert space comes with $p$ unitary operators
  $\epsilon_1,\dots,\epsilon_p$ which are odd with respect to the
  grading, which satisfy $\epsilon_i^2=-1$ and
  $\epsilon_i\epsilon_j+\epsilon_j\epsilon_i =0$ for all $i$ and $j$
  with $i\ne j$.
  An operator $T\colon H\to H$ on a $p$-multigraded Hilbert space is
  called \emph{multigraded} if it commutes with
  $\epsilon_1,\dots,\epsilon_p$. Such operators can (in addition) be
  even or odd.

  This definition can be reformulated as saying that a multigraded
  Hilbert space is a (right) module over 
  the Clifford algebra $Cl_p$, and a multigraded operator is a module
  map.

  We now define $KK_p^\Gamma(X)$ using cycles as above, with the
  additional assumption that the Hilbert space is $p$-graded, that the 
  representation $\pi$ takes values in $\pi$-multigraded even operators,
  and that the operator $F$ is an odd $p$-multigraded
  operator. Isomorphism and equivalence of these multigraded cycles is 
  defined as above, requiring that the multigradings are
  preserved throughout.
\end{definition}

This definition gives an equivariant homology theory if we restrict
to \emph{proper} maps. Moreover, it satisfies Bott
periodicity. The
period is two for the (complex) K-homology we have considered so
far. All results mentioned in this section generalize to higher
equivariant K-homology.

If $X$ is a proper $\Gamma$-CW-complex, the analytically defined representable
equivariant $K$-homology groups $RK_p^\Gamma(X)$ are canonically
isomorphic to the equivariant $K$-homology groups $K^\Gamma_p(X)$
defined by Davis and L{\"u}ck in \cite{Davis-Lueck(1998)} as described
in Section \ref{sec:homot-ther-defin}.

\subsection{The assembly map}
Here, 
  we will use the homotopy theoretic description of equivariant K-homology
  due to Davis and L{\"u}ck \cite{Davis-Lueck(1998)} described in
  Section \ref{sec:homot-ther-defin}.
  The assembly map then becomes
  particularly convenient to describe. From the present point of view,
  the main virtue
  is that they define a functor from \emph{arbitrary},
  not necessarily proper, $\Gamma$-CW-complexes to abelian groups.

  The Baum-Connes assembly map is now simply defined
  using the equivariant collapse $E(\Gamma,fin)\to *$:
  \begin{equation}\label{eq:Davis-Lueck_BC}
   \mu\colon K^\Gamma_k(E(\Gamma,fin)) \to K^\Gamma_k(*) = K_k(C^*_r\Gamma).
 \end{equation}
 If $\Gamma$ is torsion-free, then $E\Gamma=E(\Gamma,fin)$, and the
 assembly map of \eqref{eq:BC} is defined as the composition of
 \eqref{eq:Davis-Lueck_BC} with the appropriate isomorphism in
 Proposition \ref{prop:Davis_Lueck_properties}.

\subsection{Survey of KK-theory}

The analytic definition of $\Gamma$-equivariant K-homology can be extended to a 
bivariant functor on $\Gamma$-$C^*$-algebras. Here, a
$\Gamma$-$C^*$-algebra is a $C^*$-algebra $A$ with an action (by
$C^*$-algebra automorphisms) of $\Gamma$. If $X$ is a proper
$\Gamma$-space, $C_0(X)$ is such a $\Gamma$-$C^*$-algebra.

Given two $\Gamma$-$C^*$-algebras $A$ and $B$, Kasparov defines the
bivariant KK-groups $KK^\Gamma_*(A,B)$.
The most important property of this bivariant KK-theory is that it
comes with a (composition) product, the \emph{Kasparov product}. This can be
stated most conveniently as follows:

Given a discrete group $\Gamma$, we have a category $KK^\Gamma$ whose
objects are $\Gamma$-$C^*$-algebras (we restrict here to separable
$C^*$-algebras). The morphisms in this category between two
$\Gamma$-$C^*$-algebras $A$ and $B$ are called
$KK_*^\Gamma(A,B)$. They are $\integers/2$-graded abelian groups, and
the composition preserves the grading, i.e.~if $\phi\in
KK_i^\Gamma(A,B)$ and $\psi\in KK_j^\Gamma(B,C)$ then $\psi\phi\in
KK^\Gamma_{i+j}(A,C)$.

There is a functor from the category of separable
$\Gamma$-$C^*$-algebras (where morphisms are $\Gamma$-equivariant
$*$-homomorphisms) to the category $KK^\Gamma_*$ which maps an object
$A$ to $A$, and such that the image of a morphism $\phi\colon A\to B$
is contained in $KK^\Gamma_0(A,B)$.

If $X$ is a proper cocompact $\Gamma$-CW-complex then (by definition)
\begin{equation*}
KK^\Gamma_p(C_0(X),\complexs)= KK^\Gamma_{-p}(X).
\end{equation*}
Here, $\complexs$
has the trivial $\Gamma$-action.

On the other hand, for any $C^*$-algebra $A$ without a group action
(i.e.~with trivial action of hte trivial group $\{1\}$),
$KK^{\{1\}}_*(\complexs,A)= K_*(A)$.

There is a functor from $KK^\Gamma$ to $KK^{\{1\}}$, called \emph{descent},
which assigns to every $\Gamma$-$C^*$-algebra $A$ the \emph{reduced
  crossed product} $C^*_r(\Gamma,A)$. The crossed product has the
property that
$C^*_r(\Gamma,\complexs)=C^*_r\Gamma$.

\subsection{KK assembly}

We now want to give an account of the analytic definition of the
assembly map, which was the original definition. The basic idea is
that the assembly map is given by taking an index. To start with,
assume that we have an even generalized elliptic $\Gamma$-operator
$(H,\pi,F)$, representing an element in $K^\Gamma_0(X)$, where $X$ is
a proper $\Gamma$-space such that $\Gamma\backslash X$ is compact. The
index of
this operator should give us an element in $K_0(C^*_r\Gamma)$. Since
the cycle is even, $H$ split as $H= H_0\oplus H_1$, and $F=\left(
  \begin{smallmatrix}
    0 & P\\ P^* & 0
  \end{smallmatrix}\right)$ with respect to this splitting. Indeed,
now, the kernel and cokernel of $P$ are modules over
$\complexs\Gamma$, and should, in most cases, give modules over
$C^*_r\Gamma$.

If $\Gamma$ is finite, the latter is indeed the case (since
$C^*_r\Gamma=\complexs\Gamma$). Moreover, since $\Gamma\backslash X$
is compact and $\Gamma$ is finite, $X$ is compact, which implies that
$C_0(X)$ is unital. We may then assume that $\pi$ is unital (switching 
to an equivalent cycle with Hilbert space $\pi(1)H$, if
necessary). But then the axioms for a cycle imply that $F^2-1$ is
compact, i.e.~that $F$ is invertible modulo compact operators, or that
$F$ is Fredholm, which means that $\ker(P)$ and $\ker(P^*)$
are finite dimensional. Since $\Gamma$ acts on them,
$[\ker(P)]-[\ker(P^*)]$ defines an element of the representation ring
$R\Gamma=K_0(C^*_r\Gamma)$ for the finite group $\Gamma$. It remains
to show that this map respects the equivalence relation defining
$K_0^\Gamma(X)$. 

However, if $\Gamma$ is not finite, the modules $\ker(P)$ and
$\ker(P^*)$, even if they are $C^*_r\Gamma$-modules, are
in general not finitely generated projective.

To grasp the difficulty, consider Example
\ref{ex:Atiyah-operator_on_covering}. Using the description where $F$
acts on a bundle over the base space $M$ with infinite dimensional
fiber $L\tensor E$, we see that loosely speaking, the null space of
$F$ should rather
``contain'' certain copies of $l^2\Gamma$ than copies of $C^*_r\Gamma$ 
(for finite groups, ``accidentally'' these two are the
same!). However, in general $l^2\Gamma$ is not projective over
$C^*_r\Gamma$ (although it is a module over this algebra). To be
specific, assume that $M$ is a point, $E_0=\complexs$ and $E_1=0$, and
$D=0$. Here we obtain, $L^2(E_0)=l^2\Gamma$, $L^2(E_1)=0$, $F=0$, and
indeed, $\ker(P)=l^2\Gamma$.

In the situation of our example, there is a way around this problem:
Instead of twisting the operator $D$ with the flat bundle
$l^2(\Gamma)\times_\Gamma \overline M$, we twist with
$C^*_r(\Gamma)\times_\Gamma\overline M$, to obtain an operator $D'$
acting on a bundle with fiber $C^*_r\Gamma\tensor \complexs^{\dim
  E}$. This way, we replace
$l^2\Gamma$ by $C^*_r\Gamma$ throughout. Still, it is not true in
general that the kernels we get in this way are finitely generated
projective modules over $C^*_r\Gamma$. However, it is a fact that one
can always add to the new
$F'$ an appropriate compact operator such that this is the case. Then
the
obvious definition gives an element
\begin{equation*}
\ind(D')\in K_0(C^*_r\Gamma).
\end{equation*}
This is the Mishchenko-Fomenko index
of $D'$ which does not depend on the chosen compact
perturbation. Mishchenko and Fomenko give a formula for this index
extending the Atiyah-Singer index formula, compare e.g.~\cite[Section
1]{Mishchenko-Fomenko(1979)} or \cite[Section 1]{Rosenberg(1983)} and
\cite[Section 3]{Rosenberg(1987)}. .

One way to get around the difficulty in the general situation (not
necessarily studying a lifted differential operator) is to deform
$(H,\pi,F)$ to an equivalent
$(H,\pi,F')$ which is better behaved (reminiscent to the compact
perturbation above). This allows to proceeds with a rather elaborate
generalization of the Mishchenko-Fomenko
example we just considered, essentially replacing
$l^2(\Gamma)$ by $C^*_r\Gamma$ again. In this way, one defines an
index as an element of
$K_*(C^*_r\Gamma)$.

This gives a homomorphism $\mu^\Gamma\colon KK_*^\Gamma(C_0(X))\to
K_*(C^*_r\Gamma)$ for each proper $\Gamma$-CW-complex $X$ where
$\Gamma\backslash X$ is compact. This passes 
to direct limits and defines, in particular,
\begin{equation*}
  \mu_*\colon RK^\Gamma_*(E(\Gamma,fin))\to K_*(C^*_r\Gamma).
\end{equation*}

Next, we proceed with an alternative definition of the Baum-Connes map
using KK-theory
and the Kasparov product. The basic observation here is that, given
any proper $\Gamma$-CW-space $X$, there is
a specific projection $p\in C^*_r(\Gamma,C_0(X))$ (unique up to an
appropriate
type of homotopy) which gives rise to a canonical element $[L_X]\in
K_0(C^*_r(\Gamma,C_0(X)))=
KK_0(\complexs,C^*_r(\Gamma,C_0(X)))$. This defines by composition the 
homomorphism
\begin{multline*}
 KK^\Gamma_*(X)=  KK^\Gamma_*(C_0(X),\complexs)\xrightarrow{\text{descent}}
 KK_*(C^*_r(\Gamma,C_0(X)),C^*_r\Gamma)\\
   \xrightarrow{[L_X]\circ\cdot} 
  KK_*(\complexs,C^*_r\Gamma) =K_*(C^*_r\Gamma).
\end{multline*}
Again, this passes to direct limits and defines as a special case the
Baum-Connes assembly map
\begin{equation*}
\mu\colon  RK_*^\Gamma(E(\Gamma,fin)) \to K_*(C^*_r\Gamma).
\end{equation*}

\begin{remark}
  It is a non-trivial fact (due to Hambleton and Pedersen
  \cite{Hambleton-Pedersen(2001)}) that this
  assembly map coincides with the map $\mu$ of \eqref{eq:BC_general}.
\end{remark}

Almost all positive results about the Baum-Connes have been obtained
using the powerful methods of KK-theory, in particular the so called
Dirac-dual Dirac method, compare e.g.~\cite{Valette(2001)}.

\subsection{The status of the conjecture}\label{sec:status-BC-conjecture}

The Baum-Connes conjecture is known to be true for the following classes of
groups.
\begin{enumerate}
\item discrete subgroups of $SO(n,1)$ and $SU(n,1)$ \cite{Julg-Kasparov(1995)}
\item Groups with the \emph{Haagerup property}, sometimes called
  \emph{a-T-menable groups}, i.e.~which admit an
  isometric action on some affine Hilbert $H$ space which is proper,
  i.e.~such that $g_nv\xrightarrow{n\to\infty} \infty$ for every $v\in 
  H$ whenever $g_n\xrightarrow{n\to\infty} \infty$ in $G$
  \cite{Higson-Kasparov(1997)}. 
  Examples of groups with the Haagerup property are amenable groups,
  Coxeter groups, groups acting properly on trees, and groups acting
  properly on simply connected CAT(0) cubical complexes
\item One-relator groups, i.e.~groups with a presentation
  $G=\innerprod{g_1,\dots,g_n\mid r}$ with only one defining relation
  $r$ \cite{Beguni-Bettaieb-Valette(1999)}.
\item Cocompact lattices in $Sl_3(\reals)$, $Sl_3(\complexs)$ and
  $Sl_3(\rationals_p)$ ($\rationals_p$ denotes the $p$-adic numbers) \cite{Lafforgue(1999)}
\item Word hyperbolic groups in the sense of Gromov \cite{Mineyev-Yu(2001)}.
\item Artin's full braid groups $B_n$ \cite{Schick(2001a)}.
\end{enumerate}

Since we will encounter amenability later on, we recall
the definition here.
\begin{definition}\label{def:amenable}
  A finitely generated discrete group $\Gamma$ is called amenable, if
  for any given finite set of generators $S$ (where we require $1\in
  S$ and require that $s\in
  S$ implies $s^{-1}\in S$) there exists a sequence of finite subsets $X_k$
  of $\Gamma$ such that
  \begin{equation*}
    \frac{\abs{SX_k:=\{sx\mid s\in S, x\in
        X_k\}}}{\abs{X_k}}\xrightarrow{k\to\infty} 1.
  \end{equation*}
  $\abs{Y}$ denotes the number of elements of the set $Y$.
  
  An arbitrary discrete group is called amenable, if each finitely
  generated subgroup is amenable.

  Examples of amenable groups are all finite groups, all abelian,
  nilpotent and solvable groups. Moreover, the class of amenable
  groups is closed under taking subgroups, quotients,  extensions, and 
  directed unions.

  The free group on two generators is not amenable. ``Most'' examples of
  non-amenable groups do contain a non-abelian free group.
\end{definition}

There is a certain stronger variant of the Baum-Connes conjecture, the
\emph{Baum-Connes conjecture with coefficients}. It has the following
stability properties:
\begin{enumerate}
\item If a group $\Gamma$ acts on a tree such that the stabilizer of
  every edge and every vertex satisfies the Baum-Connes conjecture
  with coefficients,
  the same is true for $\Gamma$ \cite{Oyono(1998)}.
\item If a group $\Gamma$ satisfies the Baum-Connes conjecture with
  coefficients, then so does every subgroup of $\Gamma$ \cite{Oyono(1998)}
\item If we have an extension $1\to\Gamma_1\to \Gamma_2\to\Gamma_3\to
  1$, $\Gamma_3$ is torsion-free and $\Gamma_1$ as well as $\Gamma_3$
  satisfy the Baum-Connes conjecture with coefficients, then so does $\Gamma_2$.
\end{enumerate}

It should be remarked that in the above list, all groups except for
word hyperbolic groups, and cocompact subgroups of $Sl_3$ actually
satisfy the Baum-Connes
conjecture with coefficients.

The Baum-Connes assembly map $\mu$ of \eqref{eq:BC_general} is
known to be rationally injective for considerably larger classes of
groups, in particular the following.
\begin{enumerate}
\item Discrete subgroups of connected Lie groups
  \cite{Kasparov(1995)}
\item Discrete subgroups of $p$-adic groups
  \cite{Kasparov-Skandalis(1991)}
\item Bolic groups (a certain generalization of word hyperbolic
  groups) \cite{Kasparov-Skandalis(1994)}.
\item Groups which admit an amenable action on some compact space
\cite{Higson-Roe(2000)}.
\end{enumerate}

Last, it should be mentioned that recent constructions of Gromov show
that certain variants of the Baum-Connes conjecture, among them the
Baum-Connes conjecture with coefficients, and an extension called the
\emph{Baum-Connes conjecture for groupoids}, are false
\cite{Higson-Lafforgue-Skandalis(2001)}. At the
moment, no counterexample to the Baum-Connes conjecture
\ref{conj:general_BC} seems to be known. However, there are many
experts in the field who think that such a counterexample eventually
will be constructed \cite{Higson-Lafforgue-Skandalis(2001)}.

\section{Real $C^*$-algebras and K-theory}\label{sec:real-c-algebras}

\subsection{Real $C^*$-algebras}

The applications of the theory of $C^*$-algebras to geometry and
topology we present here require at some point that we work with real
$C^*$-algebras. Most of the theory is parallel to the theory of
complex $C^*$-algebras. For more details on real $C^*$-algebras and
their K-theory, including the role this plays in index theory, compare
\cite{Schroeder(1993)}.

\begin{definition}
  A unital real $C^*$-algebra is a Banach-algebra $A$ with unit over
  the real numbers,
  with an isometric involution $*\colon A\to A$, such that
  \begin{equation*}
    \abs{x}^2= \abs{x^*x}\qquad\text{and }1+x^*x\text{ is
      invertible}\quad\forall x\in A.
  \end{equation*}

  It turns out that this is equivalent to the existence of a
  $*$-isometric embedding of $A$ as a closed subalgebra into
  $\boundedops{H_\reals}$, the bounded operators on a suitable real
  Hilbert space (compare \cite{Palmer(1970)}).
\end{definition}

\begin{example}
  If $X$ is a compact topological space, then $C(X;\reals)$, the
  algebra of real valued continuous function on $X$, is a real
  $C^*$-algebra with unit (and with trivial involution).

  More generally, if $X$ comes with an involution $\tau\colon X\to X$
  (i.e. $\tau^2=\id_X$), then $C_\tau(X):=\{f\colon X\to\complexs\mid
  f(\tau x)=\overline{f(x)}\}$ is a real $C^*$-algebra with involution 
  $f^*(x)=\overline{f(\tau x)}$.

  Conversely, every commutative unital real $C^*$-algebra is
  isomorphic to some $C_\tau(X)$.

  If $X$ is only locally compact, we can produce examples of
  non-unital real $C^*$-algebras as in Example \ref{ex:loc_comp_Cstar_algebra}.
\end{example}

Essentially everything we have done for (complex) $C^*$-algebras
carries over to real $C^*$-algebras, substituting $\reals$ for
$\complexs$ throughout. In particular, the definition of the K-theory
of real $C^*$-algebras is literally the same as for complex
$C^*$-algebras (actually, the definitions make sense for even more
general topological algebras), and a short exact sequence of real
$C^*$-algebras gives rise to a long exact K-theory sequence.

The notable exception is Bott periodicity. We don't get the period
$2$, but the period $8$.

\begin{theorem}
  Assume that $A$ is a real $C^*$-algebra. Then we have a Bott
  periodicity isomorphism
  \begin{equation*}
    K_0(A)\iso K_0(S^8A).
  \end{equation*}
  This implies 
  \begin{equation*}
    K_n(A) \iso K_{n+8}(A)\qquad\text{for }n\ge 0.
  \end{equation*}
\end{theorem}

\begin{remark}
  Again, we can use Bott periodicity to define $K_n(A)$ for arbitrary
  $n\in\integers$, or we may view $K_n(A)$ as an $8$-periodic theory,
  i.e.~with $n\in\integers/8$.

  The long exact sequence of Theorem \ref{theo:long_exact} becomes a
  24-term cyclic exact sequence.
\end{remark}

The \emph{real reduced $C^*$-algebra} of a group $\Gamma$, denoted
$C^*_{\reals,r}\Gamma$, is the norm closure of $\reals\Gamma$ in the
bounded operators on $l^2\Gamma$.

\subsection{Real K-homology and Baum-Connes}

More details about the contents of this subsection can be found in
\cite[Section 2]{Rosenberg(1987)}.

A variant of the cohomology theory given by complex vector bundles is
KO-theory, which is given by real vector bundles. The homology theory
dual to this is KO-homology. If $KO$ is the spectrum of topological
KO-theory, then $KO_n(X) = \pi_n(X_+\wedge KO)$.

The homotopy theoretic definition of equivariant K-homology can be
varied easily to define equivariant KO-homology. The analytic
definition can also be adapted easily, replacing $\complexs$ by
$\reals$ throughout, using in particular real Hilbert spaces. However, 
we have to stick to $n$-multigraded cycles to define $KK^\Gamma_n(X)$, 
it is not sufficient to consider only even and odd cycles.

All the constructions and  properties translate appropriately from the 
complex to the real situation, again with the notable
exception that Bott periodicity does not give the period $2$, but
$8$. The upshot of all of this is that we get a real version of the
Baum-Connes conjecture, namely
\begin{conjecture}
  The real Baum-Connes assembly map
  \begin{equation*}
    \mu_n\colon KO^\Gamma_n(E(\Gamma,fin)) \to KO_n(C^*_{\reals,r}\Gamma),
  \end{equation*}
  is an isomorphism.
\end{conjecture}

It should be remarked that all known results about injectivity or
surjectivity of the Baum-Connes
map can be proved for the real version as well as for the
complex version, since each proof translates without
too much difficulty. Moreover, it is known that the complex version of
the
Baum-Connes conjecture for a group $\Gamma$ implies the real version
(for this abstract result, the isomorphism is needed as input, since
this is based on the use of the five-lemma at a certain point).

\bibliographystyle{plain}
\bibliography{index_theory}

\begin{thebibliography}{10}

\bibitem{Atiyah-Collected2}
Michael Atiyah.
\newblock {\em Collected works. {V}ol. 2}.
\newblock The Clarendon Press Oxford University Press, New York, 1988.
\newblock $K$-theory.

\bibitem{Baum-Connes-Higson(1994)}
Paul Baum, Alain Connes, and Nigel Higson.
\newblock Classifying space for proper actions and ${K}$-theory of group
  ${C}\sp \ast$-algebras.
\newblock In {\em $C\sp \ast$-algebras: 1943--1993 (San Antonio, TX, 1993)},
  pages 240--291. Amer. Math. Soc., Providence, RI, 1994.

\bibitem{Baum-Higson-Schick(2001)}
Paul Baum, Nigel Higson, and Thomas Schick.
\newblock Equivariant {K}-homology as via equivariant $({M,E},\phi)$-theory.
\newblock in preparation (2001).

\bibitem{Beguni-Bettaieb-Valette(1999)}
C{\'e}dric B{\'e}guin, Hela Bettaieb, and Alain Valette.
\newblock ${K}$-theory for ${C}\sp \ast$-algebras of one-relator groups.
\newblock {\em $K$-Theory}, 16(3):277--298, 1999.

\bibitem{Berline-Getzler-Vergne(1992)}
Nicole Berline, Ezra Getzler, and Mich{\`e}le Vergne.
\newblock {\em Heat kernels and {D}irac operators}.
\newblock Springer-Verlag, Berlin, 1992.

\bibitem{Berrick-Chatterji-Mislin(2001)}
A.J. Berrick, I.~Chatterji, and G.~Mislin.
\newblock From acyclic groups to the bass conjecture for amenable groups.
\newblock preprint 2001, submitted for publication.

\bibitem{Blackadar(1998)}
Bruce Blackadar.
\newblock {\em ${K}$-theory for operator algebras}.
\newblock Cambridge University Press, Cambridge, second edition, 1998.

\bibitem{Bott-Tu(1982)}
Raoul Bott and Loring~W. Tu.
\newblock {\em Differential forms in algebraic topology}.
\newblock Springer-Verlag, New York, 1982.

\bibitem{Broecker-Dieck(1995)}
Theodor Br{\"o}cker and Tammo tom Dieck.
\newblock {\em Representations of compact {L}ie groups}.
\newblock Springer-Verlag, New York, 1995.
\newblock Translated from the German manuscript, Corrected reprint of the 1985
  translation.

\bibitem{Davis-Kirk(2001)}
James~F. Davis and Paul Kirk.
\newblock {\em Lecture notes in algebraic topology}.
\newblock American Mathematical Society, Providence, RI, 2001.

\bibitem{Davis-Lueck(1998)}
James~F. Davis and Wolfgang L{\"u}ck.
\newblock Spaces over a category and assembly maps in isomorphism conjectures
  in ${K}$- and ${L}$-theory.
\newblock {\em $K$-Theory}, 15(3):201--252, 1998.

\bibitem{Hambleton-Pedersen(2001)}
Ian Hambleton and Erik Pedersen.
\newblock Identifying assembly maps in {K}- and {L}-theory.
\newblock preprint (2001), available via
  \href{http://www.math.binghamton.edu/erik/}{http://www.math.binghamton.edu/erik/}.

\bibitem{Higson-Kasparov(1997)}
Nigel Higson and Gennadi Kasparov.
\newblock Operator ${K}$-theory for groups which act properly and isometrically
  on {H}ilbert space.
\newblock {\em Electron. Res. Announc. Amer. Math. Soc.}, 3:131--142
  (electronic), 1997.

\bibitem{Higson-Lafforgue-Skandalis(2001)}
Nigel Higson, Vincent Lafforgue, and George Skandalis.
\newblock Counterexamples to the {B}aum-{C}onnes conjecture.
\newblock preprint, Penn State University, 2001, available via
  \href{http://www.math.psu.edu/higson/research.html}{http://www.math.psu.edu/higson/research.html}.

\bibitem{Higson-Roe(2000)}
Nigel Higson and John Roe.
\newblock Amenable group actions and the {N}ovikov conjecture.
\newblock {\em J. Reine Angew. Math.}, 519:143--153, 2000.

\bibitem{Higson-Roe(2001)}
Nigel Higson and John Roe.
\newblock {\em {Analytic {K}-homology.}}
\newblock {Oxford Mathematical Monographs}. Oxford University Press, Oxford,
  2001.

\bibitem{Julg-Kasparov(1995)}
Pierre Julg and Gennadi Kasparov.
\newblock Operator ${K}$-theory for the group ${\rm {s}{u}}(n,1)$.
\newblock {\em J. Reine Angew. Math.}, 463:99--152, 1995.

\bibitem{Kasparov(1995)}
G.~G. Kasparov.
\newblock ${K}$-theory, group ${C}\sp *$-algebras, and higher signatures
  (conspectus).
\newblock In {\em Novikov conjectures, index theorems and rigidity, Vol.\ 1
  (Oberwolfach, 1993)}, pages 101--146. Cambridge Univ. Press, Cambridge, 1995.

\bibitem{Kasparov-Skandalis(1991)}
G.~G. Kasparov and G.~Skandalis.
\newblock Groups acting on buildings, operator ${K}$-theory, and {N}ovikov's
  conjecture.
\newblock {\em $K$-Theory}, 4(4):303--337, 1991.

\bibitem{Kasparov-Skandalis(1994)}
Guennadi Kasparov and Georges Skandalis.
\newblock Groupes ``boliques'' et conjecture de {N}ovikov.
\newblock {\em C. R. Acad. Sci. Paris S\'er. I Math.}, 319(8):815--820, 1994.

\bibitem{Kochman(1996)}
S.~O. Kochman.
\newblock {\em Bordism, stable homotopy and {A}dams spectral sequences}.
\newblock American Mathematical Society, Providence, RI, 1996.

\bibitem{Lafforgue(1999)}
Vincent Lafforgue.
\newblock Compl\'ements \`a la d\'emonstration de la conjecture de
  {B}aum-{C}onnes pour certains groupes poss\'edant la propri\'et\'e ({T}).
\newblock {\em C. R. Acad. Sci. Paris S\'er. I Math.}, 328(3):203--208, 1999.

\bibitem{Lawson-Michelsohn(1989)}
H.~Blaine Lawson, Jr. and Marie-Louise Michelsohn.
\newblock {\em Spin geometry}.
\newblock Princeton University Press, Princeton, NJ, 1989.

\bibitem{Lueck(2001a)}
W.~L{\"u}ck.
\newblock The relation between the {B}aum-{C}onnes conjecture and the trace
  conjecture.
\newblock Preprintreihe SFB 478 --- Geometrische Strukture in der Mathematik,
  Heft 151, 2001.

\bibitem{Milnor-Stasheff(1974)}
John~W. Milnor and James~D. Stasheff.
\newblock {\em Characteristic classes}.
\newblock Princeton University Press, Princeton, N. J., 1974.
\newblock Annals of Mathematics Studies, No. 76.

\bibitem{Mineyev-Yu(2001)}
Igor Mineyev and Guoliang Yu.
\newblock The {B}aum-{C}onnes conjecture for hyperbolic groups.
\newblock preprint 2001, available via
  \href{http://www.math.usouthal.edu/~mineyev/math/}{http://www.math.usouthal.edu/~mineyev/math/}.

\bibitem{Mishchenko-Fomenko(1979)}
A.~S. Mi{\v{s}}{\v{c}}enko and A.~T. Fomenko.
\newblock The index of elliptic operators over ${C}\sp{\ast} $-algebras.
\newblock {\em Izv. Akad. Nauk SSSR Ser. Mat.}, 43(4):831--859, 967, 1979.

\bibitem{Mislin(2001)}
Guido Mislin.
\newblock Equivariant {K}-homology of the classifying space for proper actions.
\newblock Lecture notes, in preparation (2001).

\bibitem{Oyono(1998)}
Herv{\'e} Oyono-Oyono.
\newblock La conjecture de {B}aum-{C}onnes pour les groupes agissant sur les
  arbres.
\newblock {\em C. R. Acad. Sci. Paris S\'er. I Math.}, 326(7):799--804, 1998.

\bibitem{Palmer(1970)}
T.~W. Palmer.
\newblock Real ${C}\sp*$-algebras.
\newblock {\em Pacific J. Math.}, 35:195--204, 1970.

\bibitem{Rosenberg(1983)}
J.~Rosenberg.
\newblock ${C}^*$-algebras, positive scalar curvature, and the {N}ovikov
  conjecture.
\newblock {\em Publ. Math. IHES}, 58:197--212, 1983.

\bibitem{Rosenberg(1987)}
J.~Rosenberg.
\newblock ${C}^*$-algebras, positive scalar curvature, and the {N}ovikov
  conjecture iii.
\newblock {\em Topology}, 25:319--336, 1986.

\bibitem{Roy(1999)}
Ranja Roy.
\newblock The trace conjecture---a counterexample.
\newblock {\em $K$-Theory}, 17(3):209--213, 1999.

\bibitem{Schick(2001a)}
Thomas Schick.
\newblock The trace on the ${K}$-theory of group ${C}\sp *$-algebras.
\newblock {\em Duke Math. J.}, 107(1):1--14, 2001.

\bibitem{Schroeder(1993)}
Herbert Schr{\"o}der.
\newblock {\em ${K}$-theory for real ${C}\sp *$-algebras and applications}.
\newblock Longman Scientific \& Technical, Harlow, 1993.

\bibitem{Shubin(1987)}
M.~A. Shubin.
\newblock {\em Pseudodifferential operators and spectral theory}.
\newblock Springer-Verlag, Berlin, 1987.
\newblock Translated from the Russian by Stig I. Andersson.

\bibitem{Stolz(2001)}
Stephan Stolz.
\newblock Positive scalar curvature metrics on closed manifolds.
\newblock to appear in the proceedings on the School on High-Dimensional
  Manifold Topology (21 May--8 June 2001) at the abdus salam international
  centre for theoretical physics (ICTP).

\bibitem{Suslin-Wodzicki(1992)}
Andrei~A. Suslin and Mariusz Wodzicki.
\newblock Excision in algebraic ${K}$-theory.
\newblock {\em Ann. of Math. (2)}, 136(1):51--122, 1992.

\bibitem{Swan(1962)}
Richard~G. Swan.
\newblock Vector bundles and projective modules.
\newblock {\em Trans. Amer. Math. Soc.}, 105:264--277, 1962.

\bibitem{Dieck(1972)}
Tammo tom Dieck.
\newblock Orbittypen und \"aquivariante {H}omologie. {I}.
\newblock {\em Arch. Math. (Basel)}, 23:307--317, 1972.

\bibitem{Valette(2001)}
A.~Valette.
\newblock Introduction to the {B}aum-{C}onnes conjecture.
\newblock preprint 2001, to appear as ETHZ lecture note, published by
  Birkh\"auser.

\bibitem{Valette(2000)}
A.~Valette.
\newblock On the {B}aum-{C}onnes assembly map for discrete groups.
\newblock unpublished preprint, perhaps to appear as appendix to ``Introduction
  to the Baum-Connes conjecture'', to appear as ETHZ lecture note, published by
  Birkh\"auser.

\bibitem{Wegge-Olsen(1993)}
N.~E. Wegge-Olsen.
\newblock {\em ${K}$-theory and ${C}\sp *$-algebras}.
\newblock The Clarendon Press Oxford University Press, New York, 1993.
\newblock A friendly approach.

\end{thebibliography}

\end{document}